\documentclass[12pt]{amsart}
\usepackage{amssymb}
\newtheorem{theorem}{Theorem}[section]
\newtheorem{lemma}[theorem]{Lemma}
\newtheorem{proposition}[theorem]{Proposition}
\newtheorem{corollary}[theorem]{Corollary}

\theoremstyle{definition}
\newtheorem{definition}[theorem]{Definition}
\newtheorem{example}[theorem]{Example}

\theoremstyle{remark}
\newtheorem{remark}[theorem]{Remark}
\numberwithin{equation}{section}

\newcommand{\N}{\mathbb{N}}
\newcommand{\Z}{\mathbb{Z}}

\newcommand{\R}{\mathbb{R}}
\newcommand{\C}{\mathbb{C}}
\newcommand{\T}{\mathbb{T}}
\newcommand{\K}{\mathbb{K}}
\newcommand{\B}{\mathbb{B}}
\newcommand\cA{\mathcal{A}}
\newcommand\cB{\mathcal{B}}
\newcommand\cC{\mathcal{C}}
\newcommand\cE{\mathcal{E}}
\newcommand\cH{\mathcal{H}}
\newcommand\cK{\mathcal{K}}
\newcommand\cL{\mathcal{L}}
\newcommand\cO{\mathcal{O}}
\newcommand\cS{\mathcal{S}}
\newcommand\cT{\mathcal{T}}

\newcommand\Ad{\mathrm{Ad}}
\newcommand{\KE}{\K_B(\cE_E)}
\newcommand{\LE}{\mathbb{L}_B(\cE_E)}
\newcommand\inpr[2]{\langle{#1,#2}\rangle}
\newcommand\Ip{\mathrm{Ind_{\:p}\;}}
\newcommand\Ic{\mathrm{Ind_{\:cp}\;}}
\newcommand\Iw{\mathrm{Ind_{\:w}\;}}
\newcommand\In{\mathrm{Ind\;}}
\newcommand\ind{\mathrm{ind}_B^A\pi}
\newcommand\Id{\mathrm{Id}}
\newcommand\sA{A^{**}}
\newcommand\sB{B^{**}}
\newcommand\sE{E^{**}}
\newcommand\hE{\widehat{E}}
\newcommand\hsE{\widehat{\sE}}
\newcommand\Se[2]{\mathrm{Sect}(#1,#2)}
\newcommand\Mo[2]{\mathrm{Mor}(#1,#2)_0}
\newcommand\En[1]{\mathrm{End}(#1)}
\newcommand\Ho[2]{\mathrm{Hom}(#1,#2)}
\newcommand\br{\overline{\rho}}
\newcommand\bR{\overline{R}}
\newcommand\bS{\overline{S}}
\newcommand\Lr{\cL_{\rho}}

\begin{document}

\title{Inclusions of Simple $C^*$-algebras}

\author{Masaki Izumi}

\address{Department of Mathematics\\ Graduate School of Science\\
Kyoto University\\ Sakyo-ku, Kyoto 606-8502\\ Japan}
\email{izumi@kusm.kyoto-u.ac.jp}

\thanks{Work supported by Mathematical Sciences Research Institute,
JSPS, and Japanese Ministry of Education, Science, Sports and Culture}

\begin{abstract} We prove that if a conditional expectation from a
simple $C^*$-algebra onto its $C^*$-subalgebra satisfies the
Pimsner-Popa inequality, there exists a quasi-basis.
As an application, we establish the Galois correspondence for outer
actions of finite groups (and more generally finite dimensional
Kac algebras) on simple $C^*$-algebras.
We also introduce the notion of sectors for stable simple $C^*$-algebras
and purely infinite simple $C^*$-algebras in the Cuntz standard form,
and discuss several applications, including their relationship to
$K$-theory.
\end{abstract}
\maketitle

\section{Introduction} In \cite{J83}, V. F. R. Jones introduced
the notion of indices for inclusions of $\mathrm{II}_1$ factors,
which opened a completely new aspect of operator algebras involving
various other fields of mathematics and mathematical physics.
Since then, several attempts have been made to generalize the notion of
indices to broader classes of operator algebras.
In \cite{PP}, M. Pimsner and S. Popa introduced the notion of
the probabilistic index for a conditional expectation, which is the
best constant of so called the Pimsner-Popa inequality.
In \cite{K86}, H. Kosaki discussed another way to define indices
for general inclusions of factors using the theory of operator
valued weights (see also \cite{BDH}).

Though the probabilistic index works perfectly for analytic purposes
even in the case of $C^*$-algebras (see, e.g. \cite{P00}),
it is not always suitable for algebraic operations such as basic
construction in the $C^*$-case.
In \cite{W}, Y. Watatani proposed to assume existence of
a quasi-basis for a conditional expectation, a generalization of
the Pimsner-Popa basis in the von Neumann algebra case,
to analyze inclusions of $C^*$-algebras.
With a quasi-basis, Watatani successfully introduced a $C^*$-version
of basic construction \cite{W}, which is closely related to $K$-theory
of $C^*$-algebras \cite[Capter III]{W}, \cite{KajW}.
However, he also pointed out that there exists a rather annoying example
satisfying the Pimsner-Popa inequality that has no quasi-basis
\cite[Example 2.9.1]{W}.

One of the purposes of the present notes is to show that for simple
$C^*$-algebras, the Pimsner-Popa inequality always assures existence
of a quasi-basis (Theorem 3.2).
For this, we employ a second dual method, and apply M. Baillet,
Y. Denizeau, and J.-F. Havet's result \cite{BDH} in the $W^*$-case.
As far as the Pimsner-Popa inequality holds, the index of a conditional
expectation in the sense of \cite{K86} and \cite{BDH} can be defined as a
central element of the second dual.
However, it does not even belong to the multiplier algebra in general.
Indeed, our main results tell that this phenomenon occurs in the above
mentioned example.
Nevertheless, when $C^*$-algebras forming an inclusion are simple,
the index in the above sense happens to be a scalar,
and the $C^*$-basic construction works.
Using a similar method, we also give a prescription for the $C^*$-basic
construction to the non-unital case, which was not treated in \cite{W}.
As an application, we show that a minimal conditional expectation
factors through every intermediate $C^*$-subalgebra for an inclusion
of simple $C^*$-algebra with a finite index
(no matter what definition we adopt).
In particular, we establish the Galois correspondence of outer
actions of finite groups (more generally finite dimensional Kac algebras)
on simple $C^*$-algebras exactly in the same way as in the case of factors
\cite{NT}, \cite{ILP}.

As it is well-known in the case of factors, the bimodule approach
(see, for example, \cite{EKa} and references in it) to the categorical
aspect of subfactors is equivalent to the sector approach inspired
by the algebraic quantum field theory \cite{L90}.
In \cite{KajW}, T. Kajiwara and Y. Watatani developed the theory of
$C^*$-bimodules of finite indices and discuss their relationship to
$K$-theory.
We introduce here the notion of sectors for two classes of $C^*$-algebras
as an alternative approach: simple $\sigma$-unital stable
$C^*$-algebras and simple purely infinite unital $C^*$-algebras in the
Cuntz standard form.
As an application, we give a $K$-theoretical obstruction for a particular
combinatorial invariants to be realized in a simple $C^*$-algebra with
particular $K$-groups.
In the final section, we show a Kishimoto's type theorem (cf. \cite{Ki})
for conditional expectations of finite indices as anther application.

The main part of these notes is based on a series of lectures
the author gave in Tokyo Metropolitan University in December 1997,
and the author would like to thank H. Takai for giving him
the opportunity.
He is grateful to all the audiences of the lectures, and especially
H. Osaka who took notes carefully.
Some statements of this paper have already been applied in \cite{O}.
Typing of this paper was finished while the author stayed in
Mathematical Sciences Research Institute, and he would like to thank them
for their hospitality.

\section{Preliminaries}
First we summarize notations we use in this paper.
For a $C^*$-algebra $A$, $A_{sa}$ and $A_{+}$ mean the self-adjoint part
and the positive part of $A$ respectively.
We denote by $S(A)$, $P(A)$, and $Q(A)$ the set of
states on A, the set of pure states on $A$, and the set of positive linear
functional of $A$ with norm less than or equal to 1 respectively.
For $\varphi\in S(A)$, $(\pi_\varphi, H_\varphi, \Omega_\varphi)$ is the
the GNS triple of $\varphi$.
When $A$ is not unital, $M(A)$ denotes the multiplier algebra of $A$,
which is often regarded as a $C^*$-subalgebra of the second dual $\sA$.
The strict topology of $M(A)$ is defined by the family of
seminorms $\{p_a\}_{a\in A}$ where
$$p_a(x)=||ax||+||xa||. $$
The set of compact operators on $\ell^2(\N)$ is denoted by $\K$.

Throughout this section, $A\supset B$ denotes an inclusion of
$C^*$-algebras with a conditional expectation $E$ from
$A$ onto $B$.

\subsection{Indices for $C^*$-algebras}
First we introduce two numerical indices for $E$.

\begin{definition} Let $A$, $B$, and $E$ be as above.
We define
$$\Ip E=\inf\{\lambda>0;\; \Id-\frac{1}{\lambda}E \textrm{ is positive}\}.$$
$$\Ic E=\inf\{\lambda>0;\; \Id-\frac{1}{\lambda}E
\textrm{ is completely positive}\}.$$
When there are no such number $\lambda$ satisfying the first condition
(respectively the second condition), we set $\Ip E=\infty$ (respectively
$\Ic E=\infty$).
\end{definition}

$\Ip E$ was introduced by Pimsner and Popa in \cite{PP}, and is called
the {\it probabilistic index} of $E$.
Importance of the condition $\Ic E<\infty$ was pointed out by
M. Baillet, Y. Denizeau, and J.-F. Havet in \cite[Th\'eor\`em 3.5]{BDH}.
Namely, they proved that when $A$ and $B$ are von Neumann algebras,
$\Ic E<\infty$ occurs if and only if the dual operator valued weight is
bounded.

By definition, $\Ip E\leq \Ic E$ always holds.
On the other hand, it was shown by M. Frank and E. Kirchberg in
\cite[Theorem 1]{FK} that $\Ip E$ is finite if and only if $\Ic E$ is
finite.
Following their argument, we can show the following:

\begin{lemma} Let $A$, $B$, and $E$ be as above.
If $B$ has no non-zero finite dimensional representation, then $\Ip E=\Ic
E$.
\end{lemma}

\begin{proof} It suffices to show the statement when $\Ip E<\infty$, and
we assume this.
Recall that the second dual $\sB$ of $B$ is naturally regarded as a von
Neumann subalgebra of the second dual $\sA$ of $A$ \cite[Corollary
3.7.9]{Pe}.
Let $\sE$ be the normal extension of $E$ to the second dual $\sA$.
Then, $\sE$ is a conditional expectation from $\sA$ onto $\sB$ with
the same indices as $E$.
Let $z_{\sA}$ and $z_{\sB}$ be the central projections of $\sA$ and $\sB$
corresponding to the discrete parts of $\sA$ and $\sB$ respectively.
Then, thanks to \cite[Proposition 1.3]{FK}, \cite[1.1.2]{P95},
we have $z_{\sA}=z_{\sB}$.
Let $M=\sA z_{\sA}$, $N=\sB z_{\sB}$, and $F$ be the restriction of
$\sE$ to $M$.
Since the direct sum of all irreducible representations of $A$ is faithful,
it suffices to show the statement for $F$.
By assumption $B$ has no non-zero finite dimensional representation, and
$N$ is a direct sum of infinite dimensional type I factors.
Thus, there exists a system of matrix units $\{e_{ij}\}_{i,j\in \N}$ in $N$
such that $\{e_{ii}\}_{i\in \N}$ is a partition of unity.
Let $F_e$ be the restriction of $F$ to $eMe$ where $e=e_{11}$.
Then, we have
$$(M\supset N,F)\cong (eMe\otimes \B(\ell^2(\N)),
eNe\otimes \B(\ell^2(\N)),
F_e\otimes \Id).$$
Therefore, $\Ip F$ and $\Ic F$ coincide.
\end{proof}

Now we introduce an index in terms of a quasi-basis following
Watatani \cite{W}.

\begin{definition}[Unital case] Let $A$ and $B$ be unital $C^*$-algebras,
and $E$ be as above.
\begin{itemize}
\item [(1)] A quasi-basis for $E$ is a finite set
$\{(u_i,v_i)\}_{i=1}^n \subset A\times A$ such that for every $a\in A$,
the following holds:
$$a=\sum_{i=1}^n u_i E(v_i a)=\sum_{i=1}^n E(a u_i) v_i.$$
\item [(2)] When $\{(u_i,v_i)\}_{i=1}^n$ is a quasi-basis for $E$, we define
$\Iw E$ by
$$\Iw E=\sum_{i=1}^n u_iv_i.$$
When there is no quasi-basis, we write $\Iw E=\infty$.
$\Iw E$ is called the Watatani index of $E$.
\end{itemize}
\end{definition}

\begin{remark} Several remarks about the above definitions are in order.
\begin{itemize}
\item [(1)] $\Iw E$ does not depend on the choice of the quasi-basis
in the above formula, and it is a central element of $A$
\cite[Proposition 1.2.8]{W}.
\item [(2)] Once we know that there exists a quasi-basis,
we can choose one of the form $\{(w_i,w_i^*)\}_{i=1}^m$, which shows that
$\Iw E$ is a positive element \cite[Lemma 2.1.6]{W}.
\item [(3)] $\Ic E\leq ||\Iw E||$ always holds, which was first proved
in \cite{PP} in the case of factors.
\item [(4)] When, $E_1$ and $E_2$ are faithful conditional expectations
from $A$ onto $B$, then $E_1$ has a quasi-basis if and only if $E_2$ does.
In other words, the property that $E$ has a quasi-basis is actually
a property of the inclusion $A\supset B$.
\item [(5)] Though the above definition looks fine even in the non-unital
case,
Proposition 2.5 below shows that a quasi-basis exists only if both $A$ and
$B$ are unital.
It is also easy to show that a simple unitization trick does not resolve
this problem at all.
\end{itemize}
\end{remark}

Several examples of conditional expectations of Cuntz-Krieger algebras
with quasi-bases are found in \cite{I93}, \cite{I98}.

\subsection{$C^*$-basic construction}
In this subsection, we recall Watatani's notion of the $C^*$-basic
construction.
We assume that $E$ is faithful in what follows.
For general theory of Hilbert $C^*$-modules, we refer to \cite{La}.

We introduce pre-Hilbert $B$-module structure into $A$ by defining
a $B$-valued inner product
$$\inpr {a_1}{a_2}_B:=E(a_1^*a_2). $$
We denote by $\cE_E$ and $\eta_E$ the Hilbert $B$-module completion of
$A$ and the natural inclusion map from $A$ into $\cE_E$.
As usual we use notations $\KE$ and $\LE$ for the set of
``compact operators'' on $\cE_E$ and the set of bounded
$B$-module maps with adjoints respectively.
Then, the Jones projection $e_E\in \LE$ is defined by
$$e_E\eta_E(a)=\eta_E(E(a)), \quad a\in A.$$
$A$ is embedded into $\LE$ as the left multiplication operators.
Note that $\KE$ is the norm closure of the linear spans of $Ae_B A$
in $\LE$, which is called the reduced $C^*$-basic construction and is
denoted by $A_1$.

A word of caution: in general we do not have $A_1\supset A$ even when
$\Ip E<\infty$ holds.
Indeed, Watatani proved the following in \cite{W}:

\begin{proposition} Let the notations be as above.
Then,
\begin{itemize}
\item [$(1)$] Assume that  $D$ is a $C^*$-algebra with a copy of $A$
as a $C^*$-subalgebra and that there exists a projection $p\in D$
satisfying
$$pap=E(a)p,\quad a\in A,$$
such that $B\ni b \mapsto bp$ is an injective map.
Then there exists an isomorphism $\pi$ from $A_1$ into $D$ such that
$\pi$ is an identity map on $A$ and $\pi(e_B)=p$
\cite[Proposition 2.2.7]{W}.
\item [$(2)$] $E$ has a quasi-basis if and only if $\Ip E<\infty$ and
$A_1$ is unital (equivalently $\KE=\LE$) \cite[Proposition 2.1.5]{W},
\cite{R742}.
\end{itemize}
\end{proposition}

When the condition of (2) is satisfied, there exists the dual
conditional expectation $E_1$ from $A_1$ onto $A$ satisfying
$E_1(e_B)=(\Iw E)^{-1}$.
If $\{(u_i,v_i)\}_{i=1}^n$ is a quasi-basis for $E$,
$E_A$ has a quasi-basis $\{(s_i,t_i)\}_{i=1}^n$ with
$$s_i=u_ie_E(\Iw E)^{1/2},\quad t_i=(\Iw E)^{1/2}e_E v_i.$$

\subsection{Non-unital case}
In the case of non-unital $C^*$-algebras, there is no way to define indices
using a quasi-basis inside the algebra as mentioned in Remark 2.4, (5).
However, our primary needs are the dual inclusion and the dual
conditional expectation rather than the quasi-basis itself, and we still
have
a good chance to get the former even when we give up the latter.
The aim of this subsection is to supplement Watatani's $C^*$-basic
construction by using a second dual argument in the non-unital case.
We first prepare some basic facts.

We regard $\sB$ as a von Neumann subalgebra of $\sA$ as before,
and also regard $M(A)$ and $M(B)$ as $C^*$-subalgebras of $\sA$ and $\sB$
respectively.
In general, there is no inclusion relation between $M(A)$ and $M(B)$.
However, when $AB=A$ holds (which is equivalent to that some (and all)
approximate unit of $B$ is an approximate unit of $A$), we have
$M(A)\supset M(B)$ \cite[3.12.12]{Pe}.

\begin{lemma} Let $A$, $B$, and $E$ be as above satisfying $\Ip E<\infty$.
Then, the following hold:
\begin{itemize}
\item [(1)]  $A=AB$. In consequence, $M(B)$ is a unital subalgebra of
$M(A)$.
\item [(2)] The restriction of $\sE$ to $M(A)$ is a conditional
expectation from $M(A)$ onto $M(B)$.
\end{itemize}
\end{lemma}

\begin{proof} (1). Let $x\in A$.
$\Ip E<\infty$ implies $x^*x\leq (\Ip E) E(x^*x)$, and so
for $0<\theta<1/2$ there exists $u\in A$ such that $x=uE(x^*x)^\theta$
\cite[Proposition 1.4.5]{Pe}.
Therefore, the statement holds.

(2). This follows from the characterization of the self-adjoint part
$M(A)_{sa}$ in $\sA_{sa}$ as
$$M(A)_{sa}=(\widetilde{A}_{sa})^m\cap (\widetilde{A}_{sa})_m,$$
(see \cite[Theorem 3.12.9]{Pe} for the notation).
\end{proof}

Let $\pi$ be a representation of $B$ on $H$.
The induced representation $\ind$ is the representation of $A$,
(which is actually extended to $\LE$), constructed as follows \cite{R741}:
We introduce an inner product into the algebraic tensor product
$\cE_E\odot_B H_\pi$ by
$$\inpr{\xi_1\otimes_B \eta_1}{\xi_2\otimes_B \eta_2}
=\inpr{\inpr {\xi_2}{\xi_1}_B\eta_1}{\eta_2}.$$
The induced representation $\ind$ is the representation of $A$
(and also $\LE$) on the Hilbert space completion
$\cE_E \overline{\otimes}_B H_\pi$ of $\cE_E\odot_B H_\pi$
coming from the left multiplication on the first tensor component.
For $\psi\in S(A)$, we naturally have
$\ind_\psi=\pi_{\psi\cdot E}$.

\begin{lemma} Let $A$, $B$, and $E$ be as above, and $(\pi_S, H_S)$ be
the universal representation of $B$, i.e.
$$(\pi_S,H_S)=(\bigoplus_{\varphi\in S(B)}\pi_\varphi,
\bigoplus_{\varphi\in S(B)}H_\varphi ).$$
If $\Ip E$ is finite, then $\ind_S$ is quasi-equivalent to
the universal representation of $A$.
\end{lemma}

\begin{proof} It suffices to show that for every state $\varphi\in S(A)$,
$\pi_\varphi$ is contained in $\ind_S$.
Let $\psi$ be the restriction of $\varphi$ to $B$.
Since $\Ip E$ is finite, we have
$$\varphi\leq (\Ip E)\varphi \cdot E=(\Ip E)\psi\cdot E. $$
Therefore, $\pi_\varphi$ is contained in $\ind_\psi$, and we get the
result.
\end{proof}

As mentioned at the beginning of this section, the condition $\Ip E<\infty$
implies $\Ic E<\infty$, and so $\Ic \sE<\infty$.
Then, Baillet-Denizeau-Havet's result \cite[Th\'eor\`em 3.5]{BDH} shows the
following: $\cE_{\sE}$ is automatically a self-dual Hilbert $W^*$-module,
and the $W^*$-algebra $\mathbb{L}_{\sA}(\cE_{\sE})$ is the weak closure
of the linear span of $\sA e_{\sE} \sA$.
Moreover, there exists a bounded normal operator valued weight $\hsE$
from $\mathbb{L}_{\sA}(\cE_{\sE})$ to $\sA$ satisfying
$\hsE(e_{\sE})=1$.
Simple computation shows that $\hsE-\Id$ is a completely positive map.
Following H. Kosaki's idea in \cite{K86}, they defined the index of
$\sE$ by $\hsE(1)\in Z(\sA)$.
Thanks to Lemma 2.7, we may consider $\sA$ the weak closure of the image
of $\ind_S$.
Therefore, there exists a natural faithful normal representation of
$\mathbb{L}_{\sA}(\cE_{\sE})$ on
$\cE_{\sE}\overline{\otimes}_{\sB}H_S=\cE_E\overline{\otimes}_BH_S$, and
we regard $\mathbb{L}_{\sA}(\cE_{\sE})$ as a concrete von Neumann algebra
acting on $\cE_E\overline{\otimes}_BH_S$.
On the other hand, $A_1$ also naturally acts on
$\cE_E \overline{\otimes}_BH_S$.
Therefore, we may identify $e_{\sE}$ with $e_E$.
Note that $A_1$ is weakly dense in $\mathbb{L}_{\sA}(\cE_{\sE})$
under this identification.

Now, we introduce $\Iw E$ for some class of inclusions of non-unital
$C^*$-algebras.

\begin{theorem} Let $A\supset B$ be an inclusion of $C^*$-algebras,
and $E$ be a conditional expectation from $A$ onto $B$
with $\Ip E<\infty$.
Then, $A_1\supset A$ holds if and only if $\hsE(1)\in Z(M(A))$.
Moreover, if this is the case, the following hold:
\begin{itemize}
\item [$(1)$] There exists a positive bounded $A$-$A$ bimodule map
$\hE$ from $A_1$ onto $A$ satisfying
$$\hE(a_1e_Ea_2)=a_1a_2,\quad a_1, a_2\in A,$$
$$\hE(x)\geq x, \quad x\in A_{1+}.$$
For every approximate unit $\{u_\lambda\}$ of $A_1$,
$\{\hE(u_\lambda)\}$ converges to $\hsE(1)$ in the strict topology.
\item [(2)] Let $E_1$ be the conditional expectation from $A_1$ onto $A$
defined by
$$E_1(x)=\hsE(1)^{-1} \hE(x),\quad x\in A_1.$$
Then, $\K_{A_1}(\cE_{E_1})\supset A_1$ holds, i.e.
the basic construction can be repeated.
\end{itemize}
\end{theorem}

\begin{proof} Assume that $A_1\supset A$ holds.
Let $\{u_\lambda\}$ be an approximate unit of $A_1$.
Since $A_1$ is weakly dense in $\mathbb{L}_{\sA}(\cE_{\sE})$,
$\{u_\lambda\}$ strongly converges to 1 in $\mathbb{L}_{\sA}(\cE_{\sE})$,
and so $\{\hE(u_\lambda)\}$ strongly converges to $\hsE(1)$ in $A^{**}$.
However, $A_1\supset A$ implies that $\{\hE(u_\lambda a)\}$ and
$\{\hE(a u_\lambda)\}$ converge to $\hE(a)\in A$ in norm, which shows that
$\{\hE(u_\lambda)\}$ converges to $\hsE(1)\in Z(M(A))$ in the strict
topology.

Assume $\hsE(1)\in Z(M(A))$ now.
Let $a\in A$ and $\{u_\lambda\}$ be an approximate unit of $A_1$.
Then, we have
\begin{eqnarray*}||(1-u_\lambda)a||^2&=&||a^*(1-u_\lambda)^2a||
\leq ||a^*(1-u_\lambda)a||\\
&\leq &||a^*\hsE(1-u_\lambda)a||.
\end{eqnarray*}
Since $\hsE(1)\in M(A)$ and $\hsE(u_\lambda)\in A$,
$\{a^*\hsE(1-u_\lambda)a\}$ is a decreasing net in $A$
strongly converging to 0 in $\sA$.
Thus, for any $\varphi\in Q(A)$,
$$\lim_\lambda \varphi(a^*\hsE(1-u_\lambda)a)=0.$$
As $Q(A)$ is compact in weak* topology, the Dini theorem implies that
the above convergence is uniform on $Q(A)$, which means
$$\lim_\lambda ||a^*\hsE(1-u_\lambda)a||=0.$$
Thus, $a=\lim_\lambda u_\lambda a\in A_1$.

(1). We set $\hE$ to be the restriction of $\hsE$ to $A_1$.
Then, $\hE$ has the desired properties.

(2). Since $A_1$ is an ideal of $\LE$ and $A_1\supset A$, we have
$A e_E\subset A_1$.
Thanks to Lemma 2.6, (1), we have $\hsE(1)\in M(A)\subset M(A_1)$, and
so $A e_E \hsE(1)\subset A_1$.
Thus, we get
$$\K_A(\cE_{E_1})\supset Ae_E \hsE(1) e_{E_1} e_E A=Ae_E A,$$
which shows $\K_{A_1}(\cE_{E_1})\supset A_1$.
\end{proof}

\begin{remark} In the above, $\mathbb{L}_{\sA}(\cE_{\sE})
\supset \LE\supset A_1$ is identified with
$A_1^{**}\supset M(A_1)\supset  A_1$.
Indeed, $(\LE\supset A_1) \cong (M(A_1)\supset A_1)$ is standard
\cite{La}.
$A_1^{**}\cong \mathbb{L}_{\sA}(\cE_{\sE})$ follows from the fact that
$B$ is a full corner of $A_1$ and the restriction of $B$ to $\ind_S$
is quasi-equivalent to the universal representation $\pi_S$.
\end{remark}

\begin{definition} Let the notations be as above.
When $A_1\supset A$ holds,
we define the Watatani index $\Iw E$ of $E$ by $\hsE(1)\in Z(M(A))$.
We define the dual conditional expectation $E_1$ from $A_1$ onto $A$ by
$$E_1(x)=(\Iw E)^{-1}\hE(x), \quad x\in A_1.$$
$E_1$ satisfies $\Ip E_1\leq ||\Iw E||$.
When $A_1\supset A$ does not hold, we write $\Iw E=\infty$.
\end{definition}

>From Watatani's original definition, it is not so clear whether
the inclusions with finite Watatani indices are closed under natural
operations, such as taking tensor product with $\K$ or passing to
full corners.
However, thanks to Theorem 2.9, we immediately get the following because
these two operations behave very well with the second dual:

\begin{lemma} Let $A\supset B$ be an inclusion of $C^*$-algebras, and
$E$ be a conditional expectation from $A$ onto $B$ with finite
$\Iw E$.
Then,
\begin{itemize}
\item [$(1)$] Let $E\otimes\Id$ be the conditional expectation
from $A\otimes \K$ onto $B\otimes \K$.
Then, we have $\Iw (E\otimes Id)=(\Iw E)\otimes 1$.
\item [$(2)$] Let $p\in M(B)$ be a non-zero projection, such that
$pBp$ is a full corner of $B$, and $E_p$ be the restriction of
$E$ to $pAp$, which is a conditional expectation from $A$ onto $B$.
Then, we have $\Iw (E_p)=(\Iw E)p$.
\end{itemize}
\end{lemma}

\begin{proof} (1). This follows from
$(A\otimes \K)^{**}=\sA \overline{\otimes}\B(\ell^2(\N))$.
(2). This follows from the fact that when $pBp$ is a full corner of $B$,
$pAp$ is a full corner of $A$ as well, and the central supports of
$p$ in $\sB$ and $\sA$ are 1.
\end{proof}

\begin{remark} In the above situation, it is easy to show
$$M(A\otimes \K)\cap (B\otimes \K)'=(M(A)\cap B')\otimes \C.$$
When $A$ and $B$ are $\sigma$-unital $C^*$-algebras, we further have
$$M(pAp)\cap (pBp)'=(M(A)\cap B')p. $$
Indeed, since $(pBp\otimes K)$ is stable, L. G. Brown's theorem
\cite[Theorem 4.23]{B} implies that there exists an isometry
$V\in M(B\otimes \K)\subset M(A\otimes \K)$ such that $VV^*=p\otimes 1$.
Thus, $\Ad(V)$ is an isomorphism from $A\otimes \K\supset B\otimes \K$
onto $pAp\otimes \K\supset pBp\otimes \K$.
Let $x\in M(A)\cap B'$.
Then, have
$$V(x\otimes 1)V^*=VV^*(x\otimes 1)=p(x\otimes 1),$$
which shows $M(pAp)\cap (pBp)'=(M(A)\cap B')p$.
\end{remark}

We give a sufficient condition for $E$ to have finite $\Iw E$ in terms of
a quasi-basis in $M(A)$.

\begin{lemma} Let $A\supset B$ be an inclusion of non-unital
$C^*$-algebra and $E$ be a conditional expectation from $A$ onto $B$
with $\Ip E<\infty$.
We assume that there exists a sequence
$\{(u_i,v_i)\}_{i=1}^\infty\subset M(A)\times M(A)$ such that
if we set
$$T_n(a)=\sum_{i=1}^nu_iE(v_i a),\quad a\in A,$$
$\{||T_n-\Id||\}_n$ converges to zero.
Then, we have $A_1\supset A$.
In particular, if the restriction of $\sE$ to $M(A)$ has a quasi-basis,
then $A_1\supset A$ holds.
\end{lemma}

\begin{proof} Note that $M(A)$ acts on $\cE_E$ by left multiplication.
Since the norm of $A$ and $\cE_E$ are equivalent, the above condition is
equivalent to
$$\lim_{n\to \infty}||1-\sum_{i=1}^n u_i e_E v_i||=0. $$
Therefore we have
$$a_1a_2=\sum_{i=1}^\infty a_1 u_i e_E v_i a_2\in A_1,\quad a_1,a_2\in A,$$
where convergence is in norm topology.
\end{proof}

The most interesting inclusions of non-unital $C^*$-algebras
satisfying the assumption of Theorem 2.9 are those of
simple $C^*$-algebras with finite $\Ip E$, which will be discussed
in the next section.
We present two rather easy but instructive examples here.

\begin{example}
Let $X$ be a locally compact Hausdorff space with a free action of
a finite group $G$.
We set $A=C_0(X)$, $B=C_0(X/G)$, and $E$ to be the average over $G$.
Then, it is a routine work to show that $E$ has $\Iw E=\# G$, and
that the basic construction $A_1$ is the crossed product
$C_0(X)\rtimes G$ (cf. \cite[Proposition 2.8.1]{W}).
\end{example}

\begin{example} Let $M(n,\C)$ be the $n$ by $n$ matrix algebra
and $\{e(n)_{ij}\}_{i,j=1}^n$ be the canonical matrix units.
We set $\cA_n=M(n,\C)$,
$\cB_n=p_n\cA_np_n+(1-p_n)\cA_n(1-p_n)$,
and
$$\cE_n(x)=p_nxp_n+(1-p_n)x(1-p_n),$$
where $p_n=e(n)_{11}$.
Then, $\cE_n$ is a conditional expectation from $\cA_n$ onto $\cB_n$.
Simple computation shows that $\{(u_i,u_i^*)\}_{i=1}^m$ is a quasi basis
if and only if the following holds:
$$1=\sum_{i=1}^mu_ip_nu_i^*=\sum_{i=1}^m u_i(1-p_n)u_i^*.$$
It is easy to construct such elements and we get $\Iw \cE_n=2$.
Let
$$V=\sum_{i=1}^m u_ip_n\otimes e(m)_{1i}.$$
Then, we have $VV^*=1\otimes e(m)_{11}$ and $V^*V\leq p_n\otimes 1_m$.
Thus, taking the trace of the both sides, we get $n\leq m$.

Now, we set $A=\oplus_{n=1}^\infty \cA_n$,
$B=\oplus_{n=1}^\infty \cB_n$,
$E=\oplus_{n=1}^\infty \cE_n$,
where the direct sums are taken in the $C^*$-sense.
Then, it is easy to show that $M(A)=\sA$ is the von Neumann algebra
direct sum, and the same thing is true for $B$ as well.
Therefore we have $\Iw E=2$.
However, the above estimate implies that the extension of $E$ to $M(A)$
does not have a quasi-basis consisting of finitely many elements.
\end{example}

\section{The simple case}
We begin this section with recalling an example discussed in
\cite[Example 2.9.1]{W}.
Let $\lambda$ be a positive number larger than or equal to 4, and
$\{e_i\}_{i=1}^\infty$ be a sequence of projections satisfying the
Jones relation, i.e.
$$e_ie_j=e_je_i, \quad  |i-j|>1,$$
$$e_ie_{i+\pm 1}e_i=\frac{1}{\lambda}e_i.$$
Let $A\supset B$ be the $C^*$-algebra generated by
$\{1\}\cup \{e_i\}_{i=1}^\infty$ and $\{1\}\cup \{e_i\}_{i=2}^\infty$
respectively with the Markov trace $\tau$ \cite{J83}.
Then, there exists a trace preserving conditional expectation
$E$ from $A$ onto $B$.
It is shown in \cite[Example 2.9.1]{W} that although $\Ip E=\lambda$, there
exists no quasi-basis for $E$.
One big difference from the case of $\lambda<4$, where there exists a
quasi-basis, is that $A$ and $B$ are \textit{not} simple for
$\lambda\geq 4$.

The purpose of this section is to show that this type of
an annoying phenomenon never occurs when $A$ or $B$ is simple.

\subsection{Existence of Watatani indices}
The following rather easy lemma turns out to be extremely useful
for our purpose:

\begin{lemma} Let $A$ be a simple $C^*$-algebra and
$\{a_\lambda\}\subset A_+$ be a bounded increasing net converging to
$a\in Z(\sA)$ in strong topology.
Then, $a$ is a scalar and
\begin{itemize}
\item [$(1)$] if $A$ is unital, $\{a_\lambda\}$ converges in norm.
\item [$(2)$] if $A$ is non-unital, $\{a_\lambda\}$ converges
in strict topology of $M(A)$.
\end{itemize}
\end{lemma}

\begin{proof} Suppose that $a$ is not a scalar.
Then, there would exist a positive number $\epsilon$ and a non-zero
central projection $z\in Z(\sA)$ such that $za\neq 0$ and
$||az||\leq ||a||-\epsilon$.
This implies $||a_\lambda z||\leq ||a||-\epsilon$.
However, since $A$ is simple, every non-zero representation of $A$
is faithful.
Thus, we would get
$$\sup_\lambda\{||a_\lambda z||\}=\sup_\lambda\{||a_\lambda ||\}=||a||,$$
which is contradiction.
Therefore, $a$ is a scalar.
The rest of the statement follows from the same argument as in
the proof of Theorem 2.8 using the Dini theorem.
\end{proof}

\begin{theorem} Let $A$ be a simple $C^*$-algebra and $B$ be a
$C^*$-subalgebra of $A$.
If $E$ is a conditional expectation from $A$ onto $B$ with $\Ip E<\infty$,
then $\Iw E$ is finite and it is a scalar.
In particular, when $A$ is unital with $\Ip E<\infty$,
$E$ has a quasi-basis.
\end{theorem}

\begin{proof} We use the same notation as in the proof of Theorem 2.8.
Let $\{u_\lambda\}$ be an approximate unit of $A_1$.
Then, $\{\hsE(u_\lambda)\}$ is an increasing net in $A_+$ converging
to $\hsE(1)\in Z(\sA)$ in $\sA$.
Thus, Lemma 3.1 implies that $\hsE(1)$ is a scalar and
Theorem 2.8 implies that $A_1\supset A$ with
$\Iw E=\hsE(1)$.
\end{proof}

The following is a generalization of M. Rieffel's result
\cite[Therem 3.1]{R80} to the case of inclusions with finite indices.

\begin{theorem} Let $A\supset B$ be an inclusion of $C^*$-algebras
with a conditional expectation $E$ from $A$ onto $B$ with
$\Ip E<\infty$.
Then,
\begin{itemize}
\item [$(1)$] If $B$ is simple, $A$ is a finite direct sum of simple
closed two-sided ideals.
\item [$(2)$] If $A$ is simple, $B$ is a finite direct sum of simple
closed two-sided ideals.
\end{itemize}
\end{theorem}

\begin{proof} (1). First we assume that $A$ is unital.
Let $J$ be a closed two-sided ideal of $A$, and
$\{u_\lambda\}$ be an approximate unit of $J$.
Since $\{u_\lambda\}$ strongly converges to a central projection
$z$ in $Z(\sA)$ \cite[3.11.10]{Pe},
$\{E(u_\lambda)\}$ converges to a central element
$\sE(z)\in Z(\sB)$ in $\sB$.
Thus, Lemma 3.1 implies that $\sE(z)$ is a scalar and
$\{E(u_\lambda)\}$ converges in norm.
However,
$$0\leq z-u_\lambda\leq (\Ip E) \sE(z-u_\lambda)$$
implies that $\{u_\lambda\}$ converges to $z$ in norm, and so
$z\in J$ and $J=Az$.
Since the restriction of $E$ to $Z(A)$ is a state satisfying
the Pimsner-Popa inequality, $Z(A)$ is finite dimensional
\cite[Proposition 2.7.3]{W}.
Thus, $A$ is a finite direct sum of simple $C^*$-algebras.

When $A$ is non-unital, the same argument as above implies that
for every $b\in B$, $\{b^*(z-u_\lambda)^2b\}$ converges in norm.
Thanks to Lemma 2.6, (1), we have $AB=A$.
Thus, $\{u_\lambda\}$ converges to $z$ in strict topology.
The same argument as above shows that $Z(M(A))$ is finite dimensional,
and $A$ is a finite direct sum of simple $C^*$-algebras.

(2). Thanks to Theorem 3.2, we can apply (1) to $A_1\supset A$, and
so $A_1$ is a finite direct sum of simple $C^*$-algebras.
Since $A_1$ is strongly Morita equivalent to $B$, we get the result.
\end{proof}

\begin{corollary} Let $A\supset B$ be an inclusion $C^*$-algebras, and
$E$ be a conditional expectation from $A$ onto $B$ with $\Ip E<\infty$.
Assume that $B$ is simple.
Then $\Iw E$ is finite, and if moreover $A$ is unital,
$E$ has a quasi-basis.
\end{corollary}

\begin{proof} Thanks to Theorem 3.3, $A$ is a finite direct sum of simple
$C^*$-algebras.
First we assume that $A$ is unital.
Let $\{z_i\}_{i=1}^n$ be the set of minimal central projections
of $A$.
Then, $E(z_i)$ is a positive scalar, say $c_i$.
We define a conditional expectation $E_i$ from $Az_i$ onto $Bz_i$ by
$$E_i(a)=\frac{1}{c_i}E(a)z_i,\quad a\in Az_i. $$
Then, $\Ip E_i\leq c_i (\Ip E)$.
Thus, we can apply Theorem 3.2 to $E_i$ and there exist quasi-bases
$\{(u^i_k, {u^i_k}^*)\}_{k=1}^{m_k}$ for $E_i$.
$$\{(\frac{1}{\sqrt{c_i}}u_k^i,\frac{1}{\sqrt{c_i}}{u_k^i}^* )\}_{i,k}$$
is a quasi-basis for $E$ with
$$\Iw E=\sum_{i=1}^n \frac{1}{c_i}(\Iw E_i)z_i.$$

When $A$ is not unital, passing to the second dual, we get
the same formula
$$\hsE(1)=\sum_{i=1}^n \frac{1}{c_i}(\Iw E_i)z_i \in Z(M(A)),$$
where $E_i$ is defined in the same way using $z_i\in Z(M(A))$.
Thus, Theorem 2.8 shows that $\Iw E$ exists.
\end{proof}

\subsection{The stable case}
For an inclusion of non-unital $C^*$-algebras
$A\supset B$ with a conditional expectation $E$, we denote by $E^M$
the restriction of $\sE$ to $M(A)$.
Then thanks to Lemma 2.6, when $\Ip E$ is finite, $E^M$ is a conditional
expectation from $M(A)$ onto $M(B)$.
Even when $\Iw E$ is finite, Example 2.15 shows that $E^M$ does not have
a quasi-basis (consisting of finitely many elements) in general.
In this subsection, we show that for a stable and $\sigma$-unital $B$,
$\Ip E<\infty$ implies existence of a quasi-basis for $E^M$.

\begin{lemma} Let $A\supset B$ be an inclusion of $C^*$-algebras with
a conditional expectation $E$ from $A$ onto $B$.
Assume that $\Ip E$ is finite.
Then,
\begin{itemize}
\item [(1)] If $B$ is stable, so is $A$.
\item [(2)] If $B$ is $\sigma$-unital, so is $A$.
\end{itemize}
\end{lemma}
\begin{proof} Note that a $C^*$-algebra $D$ is stable if and only if
there exists a system of matrix units $\{e_{ij}\}_{i,j\in \N}$ in $M(D)$
such that
$$\sum_{i=1}^\infty e_{ii}=1$$
holds in strict topology.
Thus, (1) and (2) follow from Lemma 2.6, (1).
\end{proof}

\begin{proposition} Let $A\supset B$ be an inclusion of stable
$\sigma$-unital $C^*$-algebras with a conditional expectation
$E$ from $A$ onto $B$ with $\Ip E<\infty$.
Assume that $A$ is simple.
Then, there exist isometries $V\in M(A_1)$, $W\in M(A)$ such that
$VV^*=e_E$ and
$$e_E W=(\Iw E)^{1/2}V.$$
$\{((\Iw E)^{1/2}W^*, (\Iw E)^{1/2}W)\}$ is a quasi-basis of $E^M$
satisfying
$$E^M(WW^*)=\frac{1}{\Iw E}.$$
\end{proposition}

\begin{proof} The following is a modification of
Pimsner-Popa's argument in \cite{PP}.
We use the same notation as in the proof of Theorem 2.8.
Thanks to Theorem 3.2 and Lemma 3.5, $A_1$ is stable and
$\sigma$-unital.
Since $e_E$ is a projection in $\LE=M(A_1)$ such that $e_EA_1e_E=Be_E$
is a stable $\sigma$-unital full corner, L. G. Brown's theorem
\cite[Theorem 4.23]{B} implies that there exists an isometry $V\in M(A_1)$
satisfying $VV^*=e_E$.
We set
$$W=(\Iw E)^{-1/2}\hsE(V).$$
Note that $W$ belongs to $M(A)$ thanks to Lemma 2.6, (2) and Remark 2.9.
We claim that for every $x\in \mathbb{L}_{\sA}(\cE_{\sE})$,
$e_Ex=e_E\hsE(e_Ex)$ holds \cite{PP}.
Indeed, one can easily check it for $x\in \sA e_E\sA$.
Since $\sA e_E\sA$ is weakly dense in $\mathbb{L}_{\sA}(\cE_{\sE})$,
the claim holds.
Applying this to $V$, we get $$V=(\Iw E)^{1/2}e_E W.$$
Since $V$ is isometry, we have
\begin{equation}
1=(\Iw E) W^*e_E W.
\end{equation}
Applying $\hsE$ to the both side, we see that $W$ is an isometry.
In a similar way, using $VV^*=e_E$, we get
$$E^M(WW^*)=\frac{1}{\Iw E}.$$
Now we show the quasi-basis property.
Let $x\in M(A)$.
Then (3.1) implies
$$xe_E=(\Iw E)W^*E^M(Wx)e_E.$$
Applying $\hsE$ to the both sides, we get
$$x=(\Iw E)W^*E^M(Wx),$$
which shows that $\{((\Iw E)^{1/2}W^*, (\Iw E)^{1/2}W)\}$
is a quasi-basis for $E^M$.
\end{proof}

As a corollary we get the following:
\begin{corollary} Let $A\supset B$ be an inclusion of $\sigma$-unital
infinite dimensional $C^*$-algebras with a conditional expectation
$E$ from $A$ onto $B$.
If $A$ is simple, we have
$$\Ip E=\Ic E=\Iw E.$$
\end{corollary}

\begin{proof} We already know that if one of these quantities is
finite, so are all.
Thus, we assume that they are finite.
Thanks to Theorem 3.3, $B$ is a finite direct sum of infinite
dimensional simple $C^*$-algebras, and so $B$ does not have
a non-zero finite dimensional irreducible representation.
Therefore, Lemma 2.2 implies $\Ip E=\Ic E$.
On the other hand, in general we have $\Ip E\leq \Iw E$.

Note that $\Ic E$ and $\Iw E$ do not change if we replace
$(A,B,E)$ with
$$(A\otimes \K, B\otimes \K, E\otimes \Id).$$
Thus, we may assume that $B$ is stable.
Since $A$ is dense in $M(A)$ in the strict topology,
$\Ip E=\Ip E^M$ holds.
Therefore, Proposition 3.6 implies
$$\Iw E\leq \Ip E^M=\Ip E,$$
which finishes the proof.
\end{proof}

\subsection{The purely infinite case}
There are several properties inherited through an inclusion of a finite
index (cf. \cite{O}).
Nuclearity and exactness are typical examples of such.
We show that pure infiniteness in the sense of J. Cuntz \cite{C811}
is also inherited as an application of Theorem 3.3.
A slightly different proof using the SP-property is found in \cite{O}.
In the case of crossed product inclusions, this type of results
were proven by A. Kishimoto-A. Kumjian \cite{KK}, and
J. A. Jeong-K. Kodaka-H. Osaka \cite{JKO}.

The next lemma was announced by E. Kirchberg in 1994, which can be shown
in a similar way as in \cite[Theorem 3.2]{Ro91} using
\cite[Lemma 3.1]{Ro91}.
Let $\omega\in \beta\N\setminus \N$ be a free ultra-filter.
For a $C^*$-algebra $A$, we set
$$c_\omega(A)=\{(x_n)\in \ell^\infty(\N,A);\;
\lim_{n\to\omega} ||x_n||=0\},$$
$$A^\omega=\ell^\infty(\N,A)/c_\omega(A).$$

\begin{lemma}[Kirchberg] Let $A$ be a unital infinite dimensional
$C^*$-algebra.
Then, $A$ is purely infinite simple if and only if $A^\omega$ is simple.
\end{lemma}

\begin{theorem} Let $A\supset B$ be an inclusion of $C^*$-algebras
with a conditional expectation $E$ from $A$ onto $B$ with $\Ip E<\infty$.
If either $A$ or $B$ is simple purely infinite, the other is a finite
direct sum of simple purely infinite $C^*$-algebras.
\end{theorem}

\begin{proof} Thanks to Theorem 3.2, we may assume that $B$ is simple
purely infinite by passing to the dual inclusion if necessary.
Furthermore, Theorem 3.3 shows that we may assume that $A$ and $B$ are
unital and simple by passing to corners.
Let $E^\omega$ be the conditional expectation from $A^\omega$ onto
$B^\omega$ naturally induced by $E$.
Then, we also have $\Ip E^\omega<\infty$.
Since $B^\omega$ is simple due to Lemma 3.8, Theorem 3.3 implies that
$A^\omega$ is a finite direct sum of a simple unital $C^*$-algebras.
Suppose that $A^\omega$ is not simple.
Then, there would exist a central projection $z\in A^\omega$ with
$z\neq 0,1$.
By the usual perturbation argument, we can take a representing
sequence $(z_n)$ of $z$ such that for all $n$, $z_n$ is a projection
with $z_n\neq 0,1$.
However, since $A$ is simple, there exists a sequence
$(a_n)\in \ell^\infty(\N,A)$ such that
$||a_n||=1$, $z_na_n=a_n$, and $a_n(1-z_n)=a_n$ hold for all $n$.
Let $a=(a_n)+c_\omega(A)\in A^\omega$.
Then, $a\neq 0$ though $a$ does not commute with $z$ which is
contradiction.
Thus, $A^\omega$ is simple and $A$ is purely infinite.
\end{proof}

A unital simple purely infinite $C^*$-algebra $A$ is said to be
in the \textit{Cuntz standard from} if $[1]=0$ holds in $K_0(A)$.
In the same way as in the proof of Proposition 3.6,
(or more easily using nice $K$-theoretical
properties of unital purely infinite $C^*$-algebras shown in \cite{C811}),
we can show the following:

\begin{proposition} Let $A\supset B$ be inclusion of unital simple purely
infinite $C^*$-algebras in the Cuntz standard form with a conditional
expectation $E$ from $A$ onto $B$ with $\Ip E<\infty$.
Then, there exist isometries $V\in A_1$, $W\in A$ such that
$VV^*=e_E$ and
$$e_E W=(\Iw E)^{1/2}V.$$
$\{((\Iw E)^{1/2}W^*, (\Iw E)^{1/2}W)\}$ is a quasi-basis of $E$
satisfying
$$E(WW^*)=\frac{1}{\Iw E}.$$
\end{proposition}

\begin{proof} Since $A_1$ is also unital purely infinite simple being
strongly Morita equivalent to $B$, the only point is to show that
$[1]=[e_E]=0$ in $K_0(A_1)$.
As $A_1\supset A$ is a unital inclusion and $[1]=0$ holds in $K_0(A)$,
we have $[1]=0$ in $K_0(A_1)$.
Let $\pi$ be an isomorphism from $B$ into the corner $e_EA_1e_E$
defined by $\pi(b)=be_B$.
Then, $[e_E]=\pi_*([1_B])=0$ holds as well.
Thanks to Cuntz' observation in \cite{C811}, 1 and $e_E$ are equivalent
in $A_1$, and we get a desired isometry $V$.
The rest of the argument is the same as the proof of the Lemma 3.6.
\end{proof}

\begin{remark} Even if we drop the condition that $A$ and $B$ are in the
Cuntz standard form in the above, we still have a one-element
quasi-basis consisting of a multiple of an isometry $W$ as above.
The condition we have to give up in this case is
$$E^M(WW^*)=\frac{1}{\Iw E}.$$
Indeed, using the fact that $A_1$ is purely infinite, we can find an
isometry $V\in A_1$ satisfying $VV^*\leq e_E$, which is enough to
get such $W$.
\end{remark}

\subsection{Finite group actions}
When $A$ is a unital $C^*$-algebra and $\alpha$ is an action of
$G$ on $A$.
We say that $\alpha$ is outer if for each $g\neq e$, $\alpha_g$ is
not implemented by a unitary of $M(A)$.
In general, it is not so clear when the fixed point inclusion
$A\supset A^\alpha$ has a quasi-basis while the crossed product
inclusion $A\rtimes_\alpha G\supset A$ always does (see \cite{W}).
However, when $A$ is simple, we have the following
(cf. \cite[Proposition 2.8.6, 2.8.8]{W}):

\begin{corollary} Let $A$ be a simple $C^*$-algebra and $\alpha$ be
an action of a finite group $G$ on $A$.
Let $E$ be the condition expectation from $A$ onto the fixed point
subalgebra $A^\alpha$ given by the average over $G$.
Then,
\begin{itemize}
\item [(1)]$\Iw E$ is finite.
In particular, when $A$ is unital, $E$ has a quasi-basis.
\item [(2)] Moreover, if $\alpha$ is outer, $\Iw E=\#G$ and
the basic construction is identified with the crossed product
$A\rtimes_\alpha G$.
\end{itemize}
\end{corollary}

\begin{proof} (1). This immediately follows from Theorem 3.2.

(2). We define a projection $e\in M(A\rtimes_\alpha G)$ by
$$e=\frac{1}{\#G}\sum_{g\in G}u_g,$$
where $\{u_g\}_{g\in G}$ is the implementing unitary representation
of $\alpha$ in $M(A\rtimes_\alpha G)$.
Since $\alpha$ is outer, the crossed product $A\rtimes_\alpha G$
is simple \cite{Ki}, and so $e(A\rtimes_\alpha G)e$ is a full corner
of $A\rtimes_\alpha G$.
Thus, $A_1\supset A\supset A^\alpha$ is isomorphic to
$$e((A\otimes \B(\ell^2(G)))\rtimes_{\alpha\otimes \Ad(\rho)}G)e
\supset e(A\otimes \B(\ell^2(G)))e \supset e(A\rtimes_\alpha G)e,$$
where $\rho$ is the right regular representation of $G$.
Therefore, the statement follows from Lemma 2.11, (2).
\end{proof}

\begin{remark} When $\alpha$ is outer in the above situation,
it is straightforward to show $M(A_1)\cap A'=\C$ because we have
$A_1=A\rtimes_\alpha G$.
This also implies that the relative commutant $M(A)\cap {A^\alpha}'$
is trivial.
Indeed, for $x\in M(A)\cap {A^\alpha}'$, we can define an element
$R_x\in \LE\cap A'=M(A_1)\cap A'$ by
$$R_x\eta_E(a)=\eta_E(ax), \quad a\in A.$$
However, $R_x\in M(A_1)\cap A'=\C$ implies $x\in \C$.
\end{remark}

The above corollary assures only existence of a quasi-basis,
and says nothing about the required number of elements for
a quasi-basis in terms of the index.
Indeed, the following example shows that we have no control about it.

\begin{example} We choose a sequence $\{n_k\}_{k=1}^\infty$ of integers
larger than 1 satisfying
$$\sum_{k=1}^\infty \frac{1}{n_k}<\infty. $$
Let $A$ be the UHF algebra defined by
$$A=\bigotimes_{k=1}^\infty M(n_k,\C),$$
where $M(n,\C)$ is the $n$ by $n$ matrix algebra.
We choose a rank one projection $p_k$ from each $M(n_k, \C)$ and define
self-adjoint unitaries $u_k$ by $u_k=1-2p_k$.
We define a $G=\Z/2\Z$ action $\alpha$ by
$$\alpha=\bigotimes_{k=1}^\infty \Ad(u_k).$$
Then, $\alpha$ is an outer action of $A$.
We set $B$ to be the fixed point algebra $A^\alpha$ and define $E$ by
$$E(x)=\frac{x+\alpha(x)}{2}.$$
Then, we have $\Iw E=2$.
Since $\alpha$ is outer, $B$ is a simple $C^*$-algebra,  and
the inclusion $A\supset B$ is irreducible, i.e. $A\cap B'=\C$.
We show that the required number of elements in the quasi-basis for
$E$ depends on the sequence $\{n_k\}_{k=1}^\infty$,
and has nothing to do with the index of $E$.

Let $\tau$ be the unique normalized trace,
and $M$ be the weak closure of $A$ in the GNS representation of
$A$.
We use the same symbols $\tau$ and $\alpha$ for their normal extensions.
We denote by $||a||_2$ the $L^2$-norm of $M$ defined by $\tau$.
We set
$$v_n=\bigotimes_{k=1}^nu_k.$$
Then, simple calculation shows that for $n<m$ we have
$$||v_m-v_n||_2^2=2\big(1-\prod_{k=n+1}^m(1-\frac{2}{n_k})\big).$$
Therefore, $\{v_k\}_{k=1}^\infty$ converges to a self-adjoint unitary
operator $v\in M$ in $*$-strong topology.
Thus, the extension of $\alpha$ to $M$ is the inner automorphism
implemented by $v$.
Let $e_+$ and $e_-$ be the spectral projections of $v$ corresponding
to $1$ and $-1$ respectively.
Then, we have
$$\tau(e_-)=\frac{1}{2}\big(1-\prod_{k=1}^\infty(1-\frac{2}{n_k})\big).$$
Let $\{(u_i,u_i^*)\}_{i=1}^n$ be a quasi-basis for $E$.
Since $E$ is given by
$$E(x)=\frac{x+vxv^*}{2}=e_+xe_+  +e_-xe_-,$$
we get
$$1=\sum_{i=1}^nu_ie_+u_i^*=\sum_{i=1}^nu_ie_-u_i^*.$$
Then, the same argument as in Example 2.15 implies
$$\frac{1}{\tau(e_-)}\leq n.$$
As we can make $\tau(e_-)$ arbitrarily close to 0 by choosing the sequence
$\{n_k\}_{k=1}^\infty$, we see that there is no upper bound of
the number $n$ in terms of $\Iw E$ even for inclusions of simple
$C^*$-algebras with trivial relative commutant.
In this example, $A_1=A\rtimes_\alpha \Z/2\Z$ has two extreme traces while
$A$ has a unique trace $\tau$.
Therefore, even though the Jones projection takes value $1/2$ in the
trace extension $\tau\cdot E_1$, $e_E$ is not necessarily ``big" in
$K_0(A_1)$, which also explains this phenomenon.
\end{example}

\section{Sectors}
Theorem 3.2 enables us to perform various algebraic constructions
in the case of inclusions of simple $C^*$-algebras.
In the rest of the paper, we present such applications.
Our results seem to be new even in the case of crossed product and
fixed point inclusions for finite group actions when intermediate
$C^*$-subalgebras are considered (see Section 6).

The purpose of this section is to develop the sector theory of
$C^*$-algebras (cf. \cite{L90}).
We treat only two classes of $C^*$-algebras:
$$\cC_1:\textrm{ Simple stable $\sigma$-unital $C^*$-algebras}$$
$$\cC_2:\textrm{ Simple unital purely infinite $C^*$-algebras
in the Cuntz standard form}.$$
This is not at all a constraint because for a given inclusion of
simple $\sigma$-unital $C^*$-algebras, we can find an inclusion
consisting of $C^*$-algebras in the above classes by taking a tensor
product with $\K$ (or by passing to corners if the algebras are
purely infinite) without changing the combinatorial nature
(e.g. higher relative commutants) of the inclusions.
This is one of the non-trivial consequences of Theorem 3.2.
We mainly discuss $\cC_1$ because $\cC_2$ is treated in the same way.
Thanks to Corollary 3.7, we do not need to distinguish the three
different definitions of indices in these classes,
and we simply denote them by $\In E$.

For $A,B\in \cC_1$, we denote by $\Mo BA$ the set of
$*$-homomorphism $\rho$ from $B$ to $A$ such that the image
$\rho(B)$ has a finite index in $A$.
$\rho$ has a strictly continuous extension to $M(B)$, which is still
denoted by $\rho$ for simplicity.
We denote by $E_\rho$ the minimal expectation from $A$ onto $\rho(B)$,
which can be shown to exist in the same way as in \cite{H},
\cite[Theorem 2.12.3]{W}, using $M(A)\cap \rho(B)'$ instead of
$A\cap \rho(B)'$ (for the properties of minimal expectations, see also
\cite{KL}, \cite{L92}, \cite{KawW}).
We define the statistical dimension $d(\rho)$ of $\rho$ by
$$d(\rho)=\sqrt{\In E}.$$
$\rho$ is said to be irreducible if $M(A)\cap \rho(B)'=\C$ holds.
For two $\rho_1,\rho_2\in \Mo BA$, we define the space of intertwiners
between $\rho_1$ and $\rho_2$ by
$$(\rho_1,\rho_2)=\{V\in M(A); \; V\rho_1(x)=\rho_2(x)V,\; x\in B\}.$$
When $\rho_1$ is irreducible, $(\rho_1,\rho_2)$ has Hilbert space
structure with an inner product given by
$$V_2^*V_1=\inpr {V_1}{V_2}1,\quad V_1,V_2\in (\rho_1,\rho_2).$$

We say that $\rho_1,\rho_2\in \Mo BA$ are equivalent if there exists
a unitary $u\in M(A)$ such that $\Ad(u)\cdot \rho_1=\rho_2$.
We denote by $\Se AB$ the quotient space of $\Mo BA$ by this
equivalence relation (note that the order of $A$ and $B$ are reversed
in theses two symbols).
A  member of $\Se AB$ is said to be an $A-B$ sector.
We introduce 3 categorical operations into $\Se AB$, which is actually
similar to $KK$-operations.

The product of two sectors $[\rho_1]\in\Se AB$ and $[\rho_2]\in\Se BC$
is defined by composition:
$$[\rho_1][\rho_2]:=[\rho_1\cdot\rho_2].$$
For two $[\rho_1],[\rho_2]\in \Se AB$, the direct sum is defined
as follows:
Since $A$ is stable, there exist two isometries $S_1,S_2\in M(A)$
satisfying the $\cO_2$ relation $S_1S_1^*+S_2S_2^*=1$.
We define $\rho\in \Mo BA$ by
$$\rho(x)=S_1\rho_1(x)S_1^*+S_2\rho_2(x)S_2^*,$$
and set $[\rho_1]\oplus[\rho_2]:=[\rho]$.

\begin{lemma} Let the notations be as above.
Then, for every projection $p\in  (\rho,\rho)$, $\rho\in \Mo BA$,
there exists an isometry $V\in M(A)$ such that $VV^*=p$.
In consequence, $[\rho]$ is uniquely decomposed into
a direct sum of irreducible sectors.
\end{lemma}

\begin{proof} Let $E_p$ be the conditional expectation from
$pAp$ onto $\rho(B)p$ defined by
$$E_p(x)=\frac{1}{E_\rho(p)}E_\rho(x)p,\quad x\in pAp.$$
Then, we have $\Ip E_p<\infty$.
As $\rho(B)p$ is stable, $pAp$ is stable thanks to Lemma 3.5.
Thus, L. G. Brown's theorem \cite[Theorem 4.23]{B} implies that
there exists such $V$ as above.
Note that the space $(\rho,\rho)=M(A)\cap \rho(B)'$ is finite dimensional
\cite[Proposition 2.7.3]{W}.
Let
$$M(A)\cap \rho(B)'=\bigoplus_{i=1}^n\cA_i$$
be decomposition of $M(A)\cap \rho(B)'$ into the simple components,
and $\{e(i)_{kl}\}_{k,l=1}^{m_i}$ be a system of matrix units of
$\cA_i$.
For each $i$, we choose an isometry $V(i)\in M(A)$ satisfying
$V(i)V(i)^*=e(i)_{11}$, and set $V(i)_k=e(i)_{k1}V(i)$.
Then, $\{V(i)_k\}_{i,k}$ satisfy
$$\rho(x)=\sum_{i=1}^{n}\sum_{k=1}^{m_k}V(i)_k\rho_i(x)V(i)_k^*,$$
where $\rho_i$ is defined by
$$\rho_i(x)=V(i)^*\rho(x)V(i).$$
Since $\Ad(V(i))$ is an isomorphism from $A\supset \rho_i(B)$ onto
$e(i)_{11}Ae(i)_{11}\supset \rho(B)e(i)_{11}$, $\rho_i$ belongs to $\Mo BA$.
Therefore, we get
$$[\rho]=\bigoplus_{i=1}^n m_i[\rho_i].$$
It is easy to show that $\rho_i$ is irreducible and the decomposition
is unique.
\end{proof}

Now, we introduce conjugate sectors.
For this purpose, we construct canonical endomorphisms first \cite{L90}.

\begin{lemma} Let $A\supset B$ be an inclusion of $C^*$-algebras in
$\cC_1$ of a finite index, and $E$ be the minimal conditional
expectation from $A$ onto $B$.
Then, there exists an isomorphism $\gamma_1$ from $A_1$ onto $B$ such that
if $\gamma$ is the restriction of $\gamma_1$ to $A$, there exist two
isometries $V\in (\Id_{A_1},\gamma_1)$ and $W\in (\Id_A,\gamma)$
satisfying
$$V^*W=\gamma_1(V)^*W=\frac{1}{\sqrt{\In E}},$$
$$E_1(x)=W^*\gamma_1(x)W,\quad x\in A_1.$$
$$E(a)=\gamma_1(V^*)\gamma(a)\gamma_1(V),\quad a\in A.$$
\end{lemma}

\begin{proof} Let $V$ and $W$ be the isometries defined in
Proposition 3.6.
Then, we have $V^*W=(\In E)^{-1/2}$ and $E_1^M(V)=(\In E)^{-1/2}W$.
Since $e_EA_1e_E=Be_E$ and $VV^*=e_E$ hold, there exists an
surjective isomorphism from $\gamma_1$ satisfying
$VxV^*=\gamma_1(x)e_E$ for $x\in A_1$, which implies
\begin{equation}
Vx=\gamma_1(x)V.
\end{equation}
Therefore, when $\gamma_1$ is regarded as a map from $A_1$ into
itself, $V\in (\Id_A,\gamma_1)$.
Applying $E_1^M$ to the both sides of (4.1) for $x\in A$, we get
$W\in (\Id, \gamma)$.
Note that $\gamma_1$ is given by
$$\gamma_1(x)=(\In E)E_1^M(VxV^*),\quad x\in A_1.$$
Thus,
\begin{eqnarray*} \gamma_1(V^*)W&=&(\In E)E_1^M(VV^*V^*)W
=\sqrt{\In E}E_1^M(VV^*)\\
&=&\frac{1}{\sqrt{\In E}}.
\end{eqnarray*}
For $a_1,a_2\in A$, we have
$$W^*\gamma_1(a_1e_Ea_2)W=a_1W^*\gamma_1(VV^*)Wa_2
=\frac{1}{\In E}a_1a_2=E_1(a_1e_Ea_2),$$
which shows $E_1(x)=W^*\gamma_1(x)W$.
As Proposition 3.6 implies $A=BWW^*B$, to prove $E(a)=\gamma_1(V^*aV)$,
it suffices to show it for $a=b_1WW^*b_2$, $b_1,b_2\in B$.
Since $\gamma_1$ is a map onto $B$, there exist $x_1,x_2\in A_1$
such that $\gamma_1(x_i)=b_i$, $i=1,2$.
Thus, we get
\begin{eqnarray*}
\gamma_1(V^*b_1WW^*b_2V)&=&\gamma_1(x_1V^*WW^*Vx_2)
=\frac{1}{\In E}\gamma_1(x_1x_2)\\
&=&E(b_1WW^*b_2),
\end{eqnarray*}
which finishes the proof.
\end{proof}

\begin{remark} It was pointed out in \cite[Theorem 5.12]{KajW}
that in general, downward basic construction does not exist even
for an inclusion of simple $C^*$-algebras of a finite index.
However, the above lemma shows that when $B$ is stable,
it always exists.
Indeed, we define a conditional expectation $E_{-1}$ from
$B$ onto $\gamma(A)$ by $\gamma\cdot E_1\cdot \gamma_1^{-1}$ and set
$e_{E_{-1}}=WW^*$.
Then, we have $Be_{E_{-1}}B=A$ and
$$e_{E_{-1}}be_{E_{-1}}=WW^*\gamma_1(\gamma_1^{-1}(b))WW^*
=WE_1(\gamma_1^{-1}(b))W^*=E_{-1}(b)e_{E_{-1}},\quad b\in B.$$
Therefore, we get the downward tower
$$A\supset B\supset \gamma(A)\supset \gamma(B)\supset \cdots.$$
\end{remark}

$\gamma$ is called the \textit{canonical endomorphism} for $A\supset B$.
In the above argument, the only choice made was that of $V$, and
$V$ could be replaced with $uV$, where $u$ is a unitary in $M(B)$.
Thus, $\gamma$ is uniquely determined by the inclusion as a
$B-A$ sector.
When $\gamma$ is regarded as a map from $B$ to itself,
it induces a map on the $K$-groups of $B$, which is nothing  but
the transfer map discussed in \cite{W} and \cite{KajW}.

\begin{lemma} For $A, B\in \cC_1$ and $\rho\in \Mo BA$ there exists
$\br\in \Mo AB$, unique up to equivalence, such that the following
holds: there exist isometries
$R_\rho\in (\Id_B, \br\cdot\rho)$ and
$\bR_\rho\in (\Id_A, \rho\cdot \br)$
such that
$$\bR^*\rho(R)=R_\rho^*\br(\bR_\rho)=\frac{1}{d(\rho)}.$$
Moreover, $d(\rho)=d(\br)$ holds.
\end{lemma}

\begin{proof} We apply the previous lemma to $A\supset \rho(B)$
and set $\br=\rho^{-1}\cdot \gamma$,
$R_\rho=\rho^{-1}(\gamma_1(V))$,
and $\bR_\rho=W$.
Then, these have the desired properties.
Note that the conditions in the statement imply that
$$B\supset \br(A)\supset \br\cdot\rho(B)$$
is the basic construction
(regardless the particular construction we made).
Therefore, we get $d(\br)=d(\rho)$.

Now, we show uniqueness.
Assume that $(\sigma, S,\bS)$ satisfies the same property
for $(\br, R_\rho,\bR_\rho)$ in the statement.
Then, $\sigma(\bR_\rho^*)S\in (\br,\sigma)$.
We show that $\sigma(\bR_\rho^*)S$ is a non-zero multiple of
a unitary, which will finish the proof.
We set $E(x)=\rho(S^*\sigma(x)S)$, $x\in A$.
Then, $E$ is a conditional expectation from $A$ onto $\rho(A)$
with a quasi-basis $\{(d(\rho)S^*,d(\rho)S)\}$.
Since $\In E=\In E_\rho$ and $E_\rho$ is the minimal expectation,
we get $E=E_\rho$.
Thus,
$$\rho(S^*\sigma(\bR_\rho \bR_\rho^*)S)=E_\rho(\bR_\rho \bR_\rho^*)
=\rho(R_\rho^* \br(\bR_\rho \bR_\rho^*)R_\rho)=\frac{1}{d(\rho)^2},$$
which shows that $d(\rho)\sigma(\bR_\rho^*)S$ is an isometry.
In a similar way, we get $E_\sigma(x)=\sigma(\bS^*\rho(x)\bS)$,
$x\in B$.
Therefore,
$$E_\sigma(d(\rho)^2\sigma(\bR_\rho^*)SS^*\sigma(\bR_\rho))
=d(\rho)^2\sigma(\bR_\rho^*)E_\sigma(SS^*)\sigma(\bR_\rho)
=1.$$
Since $E_\sigma$ is faithful, $d(\rho)\sigma(\bR_\rho^*)S$ is a
unitary.
\end{proof}

The sector $[\br]$ is called the \textit{conjugate sector} of $[\rho]$.

With isometries $R_\rho$, $\bR_\rho$, it is possible to show that
basically, every combinatorial result known for factors holds also
for inclusions of $C^*$-algebras in $\cC_1$, such as the Frobenius
reciprocity discussed in \cite{I98} (c.f. \cite{KajW}).
Moreover, Remark 2.11 shows that taking tensor product with
$\K$ does not change higher relative commutants, which allows us
to talk about principal graphs, standard invariants and so on.

\begin{definition} Let $A\supset B$ be an inclusion of simple
$\sigma$-unital $C^*$-algebras with a finite index.
We denote by $\iota$ the inclusion map from $B\otimes \K$ into
$A\otimes \K$.
We call the set of sectors generated by $[\iota]$ and
$[\overline{\iota}]$ the \textit{sectors associated with $A\supset B$}.
When it only contains finitely many irreducibles, $A\supset B$ is said to
be of \textit{finite depth}.
\end{definition}

As already pointed out in \cite{W}, \cite{KajW},
the most significant part of sectors in the $C^*$-case is that
they give rise to $KK$-elements, which obey particular fusion rules if
one starts from an inclusion with a particular combinatorial nature.

\begin{remark}
All the arguments in this section also work for inclusions of
$C^*$-algebras in the class $\cC_2$.
The only difference in the proof is to use the following
observation, which was already used in Proposition 3.9:
When $A\supset B$ is an inclusion of $C^*$-algebras in $\cC_2$ of
a finite index, $[p]=0$ holds in $K_0(A)$ for every projection
$p\in A\cap B'$.
This follows from the fact that $p$ is the image of the
unit of $B$ via the inclusion map of $pAp\supset Bp$.
\end{remark}

As a corollary of Lemma 4.2 (for the $\cC_2$ case), we get
the following, which is announced in \cite{P95} too without a proof:

\begin{corollary} Let $M\supset N$ be an inclusion of type III factors
with separable preduals.
If a conditional expectation $E$ from $M$ onto $N$ satisfies
$\Ip E<\infty$, $E$ is automatically normal.
\end{corollary}

\begin{proof} Note that every type III factor belongs to $\cC_2$.
Lemma 4.2 implies that there exist an endomorphism $\gamma$ of
$M$ and  an isometry $R\in M$ such that we have
$$E(x)=R^*\gamma(x)R,\quad x\in M.$$
Thanks to Takesaki's theorem \cite{T}, every endomorphism of
a von Neumann algebra with a separable predual is automatically
normal, and we get the result.
\end{proof}

\begin{example} Let $A\supset B$ be an inclusion of simple
$\sigma$-unital $C^*$-algebras with the principal graph $A_4$
(or in other words, of an index $(3+\sqrt{5})/2$).
Then, $[\iota][\overline{\iota}]$ is decomposed into two sectors
$[\Id]$ and $[\sigma]$, and $[\sigma]$ satisfies the following fusion rule
\cite{I91}:
$$[\sigma]^2=[\Id]\oplus [\sigma].$$
This implies that there exist two elements $x_i\in \En{K_i(A)}$, $i=0,1$
such that $x_i^2=x_i+1$ holds.
More specifically, we assume that $A$ is stably isomorphic to the Cuntz
algebra $\cO_{n+1}$ \cite{C77}.
Then, we have $K_0(A)\cong \Z/n\Z$, $K_1(A)=\{0\}$.
Therefore the Diophantine equation
$$x^2\equiv x+1, \mod n$$
should have a solution.
Here is the list of small integers having this property:
$$1,\;5,\;11,\;19,\;29,\;31,\;41, \cdots.$$

Now we construct examples realizing some of the above groups.
In \cite[Example 3.1]{I93}, an endomorphism $\rho\in \En{\cO_2}$
whose image has an index
$(3+\sqrt{5})/2$ is constructed by the following formula:
$$\rho(S_1)=\frac{1}{d}S_1+\frac{1}{\sqrt{d}}S_2S_2,$$
$$\rho(S_2)=(\frac{1}{d}S_1-\frac{1}{d}S_2S_2)S_2^*+S_2S_1S_1^*.$$
where $d=(1+\sqrt{5})/2$.
$\rho$ and $E_\rho$ commute with the following $\T=\R/2\pi \Z$-action:
$$\lambda_t(S_1)=e^{2it}S_1,\quad \lambda_t(S_2)=e^{it}S_2.$$
Let $G_n$ be a cyclic subgroup of $\T$ of order $n$.
Then, the restriction of $E_\rho$ to the fixed point algebra
$\cO_2^{G_n}$ is a conditional expectation onto $\rho(\cO_2^{G_n})$.
Since $E_\rho(S_1S_1^*)=1/d^2$ holds and $S_1S_1^*\in \cO_2^{G_n}$,
the index of the restriction is the same as $\In E=(3+\sqrt{5})/2$.
It is a routine work to show that $\cO_2^{G_n}$ is stably isomorphic
to the Cuntz-Krieger algebra $\cO_{D^n}$ \cite{CK}, where
$$D=\left[\begin{array}{cc}1&1\\1&0
\end{array}\right].$$
Thus, thanks to \cite{C812}, the $K$-groups of $\cO_2^{G_n}$
are given by
$$K_0(\cO_2^{G_n})=\Z^2/(1-D^n)\Z^2,$$
$$K_1(\cO_2^{G_n})=\{0\}.$$
For example,
$$K_0(\cO_2^{G_2})=\{0\},\quad K_0(\cO_2^{G_3})=\Z/2\Z\oplus\Z/2\Z,
\quad K_0(\cO_2^{G_4})=\Z/5\Z,$$
$$K_0(\cO_2^{G_5})=\Z/11\Z,\quad K_0(\cO_2^{G_6})=\Z/4\Z\oplus\Z/4\Z,
\quad K_0(\cO_2^{G_7})=\Z/29\Z.$$

It is possible to construct various inclusions of Cuntz-Krieger algebras
for each paragroup in a similar way using \cite{I98}.
For example, the principal graph $E_8$ can be realized in
the Cuntz-Krieger algebra (in fact, it is stably isomorphic to
the Cuntz algebra $\cO_{32}$ thanks to \cite{Ro95})
with $K_0$ isomorphic to $\Z/31\Z$.
Since descendant sectors of the $E_8$ inclusion contains a sector
with the principal graph $A_4$ \cite{I91}, we see that $\Z/31\Z$ can
be also realized in the above list.
\end{example}

\section{Normal extensions of conditional expectations}
There seems to be a common belief, at least among some specialists,
saying that ``when a conditional expectation leaves a state invariant,
it should have a normal extension to the weak closure in the GNS
representation of the state''.
Probably, it comes from an analogy to the case of group actions,
or depth 2 inclusions.
In this section, we provide a counter example for this statement as
an application of sector theory developed in the last section.

The following shows that there is no real distinction between induction
and reduction when an inclusion has a finite index.

\begin{lemma} Let $A$, $B$ be $C^*$-algebras in $\cC_1$
(or $\cC_2$), and $\sigma$, $\pi$ be representations of $A$ and $B$
respectively.
\begin{itemize}
\item [$(1)$] Let $\rho\in \Mo BA$, then the map
$$\Ho {\sigma\cdot \rho}{\pi}\ni X\mapsto X\sigma(\bR_\rho)\in
\Ho{\sigma}{\pi\cdot \br}$$ is a vector space isomorphism.
\item [$(2)$] If $A\supset B$ has a finite index, the following
holds:
$$\sigma|_B \cong \mathrm{ind}_A^{A_1}\sigma\cdot \gamma_1^{-1},$$
$$\ind\cong \pi\cdot\gamma,$$
where $\gamma$ and $\gamma_1$ are as in Lemma 4.2.
\end{itemize}
\end{lemma}

\begin{proof}
(1) It is straightforward to show that the inverse map is given by
$$ \Ho{\sigma}{\pi\cdot \br}\ni Y \mapsto d(\rho)\pi(R_\rho^*)Y\in
\Ho {\sigma\cdot \rho}{\pi},$$
which shows the statement.

(2). Since $A_1=V^*VA_1\subset V^*e_EAe_EA=V^*A$,
there exists an isometry $U$ from $H_\sigma$ onto
$\cE_{E_1}\overline{\otimes}_AH_\sigma$ such that
$$U\sigma(a)\xi=\sqrt{\In E}V^*a\odot_A \xi,
\quad a\in A,\;\xi\in H_\sigma,$$
which shows the first statement.
The second one can be shown in the same way using $A=W^*B$.
\end{proof}

\begin{lemma}
Let $A\supset B$ be an inclusion of $C^*$-algebras
with a conditional expectation $E$ from $A$ onto $B$ satisfying
$\Ip E<\infty$.
For every pure state $\psi\in P(B)$, there exist pure state extensions
$\{\varphi_i\}_{i=1}^n\subset P(A)$ with mutually orthogonal supports
in $\sA$ and positive numbers $\{c_i\}_{i=1}^n$ such that
$$\psi\cdot E=\sum_{i=1}^n c_i\varphi_i.$$
Moreover, $(\Ip E)^{-1}\leq c_i$ and $n\leq \Ip E$ hold.
In particular, we have
$$\ind_\psi\cong \bigoplus_{i=1}^n\pi_{\varphi_i}.$$
\end{lemma}

\begin{proof} We denote by $S_\psi(A)$ and $P_\psi(A)$ the set of state
extensions of $\psi$ to $A$, and that of the pure state extensions of
$\psi$ to $A$ respectively.
$P_\psi(A)$ is the set of extreme points of $S_\psi(A)$, and
$S_\psi(A)$ is the weak$*$ closure of the convex hull of $P_\psi(A)$.
For a state $\varphi\in S(A)$ and a representation $\pi$
of $A$, we denote by $s(\varphi)$ and $c(\pi)$ the support of $\varphi$ in
$\sA$ and the central support of $\pi$ in $Z(\sA)$.
We use the same symbols $\varphi$ and $\pi$ for their normal extensions to
$\sA$ for simplicity.

Let $\varphi\in P_\psi(A)$.
Then, we have
$$\psi\cdot E=\varphi\cdot E\geq \frac{1}{\Ip E}\varphi.$$
In particular,
\begin{equation}
\psi\cdot E(c(\pi_\varphi))\geq
\psi\cdot E(s(\varphi)) \geq \frac{1}{\Ip E},
\end{equation}
which proves the second part of the statement.

(5.1) implies that there exist finitely many mutually disjoint
irreducible representations $\pi_1,\pi_2,\cdots, \pi_m$ of $A$
such that whenever $\omega\in P_\psi(A)$, $\pi_\omega$ is equivalent
to  $\pi_j$ for some $j$.
We may assume $\psi\cdot E(c(\pi_j))\neq 0$ for all $j$.
Let $h_j$ be the density matrix for $\psi\cdot E(c(\pi_j)\cdot)$,
that is, $h_j$ is a trace class operator in $\B(H_{\pi_j})$ satisfying
$$\psi\cdot E(c(\pi_j)x)=\mathrm{Tr}(h_j\pi_j(x)),
\quad x\in A.$$
Let $s(h_j)$ be the support of $h_j$.
Since $\psi$ is pure, any state $\omega\in S(A)$ satisfying
$cs(h_j)\geq s(\omega)$ for some number $c>0$ is in $S_\psi(A)$.
Thus, (5.1) implies that $h_j$ is a finite rank operator.
We define a finite dimensional $C^*$-algebra $D$ by
$$D=\bigoplus_{j=1}^m \B(s(h_j)H_{\pi_j}).$$
Then, (5.1) also implies that there exists a continuous convex map
$\theta: S(D)\longrightarrow S_\psi(A)$ given by
$$\theta(\phi)(x)=\phi(\bigoplus_{j=1}^m s(h_j)\pi_j(x)s(h_j)),\quad
\phi\in S(D),\;x\in A.$$
Since $S(D)$ is compact, the image of $\theta$ is compact
in weak$*$ topology.
Moreover, as the image of $\theta$ contains $P_\psi(A)$,
$\theta$ is a surjection, which implies
$$\psi\cdot E(x)=\sum_{j=1}^m\psi\cdot E(c(\pi_j)x)
=\sum_{j=1}^m\mathrm{Tr}(h_j\pi_j(x)),\quad x\in A.$$
Further decomposing $h_j$, we get the decomposition of
the desired form.
\end{proof}

The following theorem tells that the depth 2 case is a rather
fortunate exception.

\begin{theorem} Let $A\supset B$ be an inclusion of $C^*$-algebras in
either $\cC_1$ or $\cC_2$, and $E$ be a conditional expectation from
$A$ onto $B$ with a finite index.
\begin{itemize}
\item [(1)] For a state $\psi\in S(B)$, there exists a normal
extension of $E$ to $\pi_{\psi\cdot E}(A)''$ if and only if
$\pi_\psi\cdot \gamma$ is quasi-equivalent to $\pi_\psi\cdot \gamma^2$,
where $\gamma$ is the canonical endomorphism for $A\supset B$.
\item [(2)] Let $\rho\in \En{\cO_2}$ be as in Example 4.8.
Then, there exists an $E_\rho$-invariant pure state $\varphi$ such that
$E_\rho$ does not have a normal extension to $\pi_\varphi(\cO_2)''$.
\end{itemize}
\end{theorem}

\begin{proof} (1). Since we have
$$\gamma(x)=(\In E)E(WxW^*),\quad x\in A,$$
$$E(x)=\gamma_1(V)^*\gamma(x)\gamma_1(V), \quad x\in A,$$
$E$ has a normal extension if and only if $\gamma$ does,
which is further equivalent to that $\pi_{\psi\cdot E}$ is
quasi-equivalent to $\pi_{\psi\cdot E}\cdot \gamma$.
On the other hand, Lemma 5.1, (2) implies
$$\pi_{\varphi\cdot E}\cong \ind_\psi\cong\pi_\psi\cdot \gamma,$$
which shows the statement.

(2). Note that $[\rho]$ is a self-conjugate sector obeying
the fusion rule
\begin{equation}
[\rho]^2=[\Id]\oplus [\rho].
\end{equation}
We claim that there exist two disjoint irreducible
representations $\pi_1$, $\pi_2$ of $A$ such that
$\pi_1\cdot\rho\cong \pi_2$ and
$\pi_2 \cong \pi_1\oplus\pi_2$.
We choose an arbitrary irreducible representation $\pi$.
Since $\In E_\rho=(3+\sqrt{5})/2<3$, Lemma 5.2 implies
that $\pi\cdot \rho$ is decomposed into at most two irreducible
representations, and there  are two cases: (i) $\pi\cdot \rho$ is
irreducible. (ii) There exist two irreducible representations
$\sigma_1$, $\sigma_2$ such that
$\pi\cdot\rho\cong \sigma_1\oplus \sigma_2$.
We first consider the case (i).
We set $\pi_1=\pi$, $\pi_2=\pi\cdot \rho$.
Then, (5.2) implies
$$\pi_2\cdot\rho=\pi\cdot \rho^2\cong \pi\oplus\pi\cdot \rho
\cong \pi_1\oplus\pi_2,$$
which also shows $\pi_1$ and $\pi_2$ are disjoint.
We consider the case (ii) now.
Lemma 5.1, (1) implies that $\sigma_1\cdot \rho$ and
$\sigma_2\cdot \rho$ contain $\pi$.
On the other hand, as before (5.2) implies
$\pi\cdot\rho^2\cong \pi\oplus \sigma_1\oplus\sigma_2$.
This means that either $\sigma_1$ or $\sigma_2$ is equivalent to
$\pi$ and we may assume $\sigma_2\cong \pi$.
Therefore, we get
$\pi\cdot\rho\cong \sigma_1\oplus \pi$ and
$\sigma_1\cdot \rho\cong \pi$.
$\pi_1:=\sigma_1$ and $\pi_2:=\pi$ satisfy the claim.

Let $A=\cO_2$ and $B:=\rho(A)$.
We take a pure state $\psi\in P(B)$ satisfying
$\pi_\psi\cong \pi_1\cdot\rho^{-1}$, and set
$\varphi=\psi\cdot E$.
Then, since $\gamma=\rho^2$, we have
$$\pi_\varphi\cong \pi_\psi\cdot \rho^2\cong \pi_1\cdot \rho
\cong \pi_2,$$
which shows that $\varphi$ is pure.
We also have
$$\pi_\psi\cdot \gamma^2\cong \pi_2\cdot\rho^2\cong
\pi_1\oplus \pi_2\oplus \pi_2.$$
This shows that $\pi_\varphi$ and $\pi_\varphi\cdot \gamma$
are not quasi-equivalent, and $E_\rho$ does not have a normal
extension to $\pi_\varphi(A)''=\B(H_\varphi)$.
\end{proof}

\section{Intermediate $C^*$-subalgebras}
In this section, following the strategy developed in \cite{ILP},
we investigate the structure of intermediate $C^*$-subalgebras of
an inclusion of simple $C^*$-algebras of a finite index.

\begin{proposition} Let $A\supset B$ be an inclusion of simple
$\sigma$-unital $C^*$-algebras, and $E$ be a conditional expectation
from $A$ onto $B$ of a finite index.
Then, for every intermediate $C^*$-algebra $A\supset C\supset B$,
there exists a conditional expectation $F$ from $A$ onto $C$
of a finite index.
\end{proposition}

\begin{proof} Let $e$ be a minimal projection of $\K$.
If there exists a conditional expectation $F$ from
$A\otimes\K$ onto $C\otimes \K$ of a finite index, the restriction
of $F$ to $(1\otimes e)(A\otimes \K)(1\otimes e)$ gives a desired
expectation.
Therefore, we may assume that $A$, $B$, and $C$ are stable from the
beginning.
Thanks to Proposition 3.6 and the argument in the proof of Corollary 3.4,
there exists a quasi-basis $\{(u_i,u_i^*)\}_{i=1}^n$ for the
restriction of $E^M$ to $M(C)$.
We define a completely positive map $F$ from $A$ to $C$ by
$$F(x)=(\Iw E|_C)^{-1}\sum_{i,j=1}^nu_iE(u_i^*xu_j)u_j^*,$$
which satisfies the following for $a\in A_+$:
\begin{eqnarray*}
F(a)&\geq&
(\Iw E|_C)^{-1}\sum_{i,j=1}^nu_iE(u_i^*xu_j)u_j^*(\Iw E|_C)^{-1}\\
&\geq&
\frac{1}{\In E}(\Iw E|_C)^{-1}\sum_{i,j=1}^nu_iu_i^*au_ju_j^*
(\Iw E|_C)^{-1}\\
&=&\frac{1}{\In E}a.
\end{eqnarray*}
We show that $F$ is a conditional expectation from $A$ onto $C$.
Let $\{(v_k,v_k^*)\}_{k=1}^m$ be another quasi-basis for $E^M|_{M(C)}$.
Then, we have
\begin{eqnarray*}
F(x)&=&\sum_{k=1}^m\sum_{i,j=1}^nv_kE^M(v_k^*u_i)E(u_i^*xu_j)u_j^*\\
&=&\sum_{k=1}^m\sum_{i,j=1}^nv_kE(E^M(v_k^*u_i)u_i^*xu_j)u_j^*\\
&=&\sum_{k=1}^m\sum_{j=1}^nv_kE(v_k^*xu_j)u_j^*.
\end{eqnarray*}
Let $u$ be a unitary in $M(C)$.
Since $\{(u^*u_i,u_i^*u)\}_{i=1}^n$ is also a quasi basis for
$E^M|_{M(C)}$, we get $F(ux)=uF(x)$, which shows that $F$ is a
$C$-$C$ bimodule map.
$F^{**}(1)=1$ implies that $F$ has norm 1.
Since $F$ is identity on $C$, $F$ is a conditional expectation.
\end{proof}

To give the structure theorem of intermediate $C^*$-subalgebras,
we briefly explain the ``crossed product decomposition" of $A$
introduced in \cite{ILP}.
As mentioned before, the Frobenius reciprocity for sectors can be
established in the same way as in \cite[Section 2]{I98}, and we use it
freely in what follows.
Let $A\supset B$ be an inclusion of $C^*$-algebras in $\cC_1$, and
$E$ be a conditional expectation from $A$ onto $B$ of a finite index.
We assume that $A\supset B$ is irreducible, that is, $M(A)\cap B'$
is trivial.
We denote by $\iota$ the inclusion map from $B$ into $A$.
Let
$$[\overline{\iota}][\iota]=\bigoplus_{r\in \cS}n_r[\rho_r],$$
be the irreducible decomposition.
We arrange the index set $\cS$ so that we have $\rho_e=\Id_B$
and $\overline{[\rho_r]}=[\rho_{\overline{r}}]$.
For simplicity, we write $d(r)=d(\rho_r)$.
For $r\in \cS$, we set
$$\cH_r=(\iota,\iota\cdot \rho_r).$$
Then, the Frobenius reciprocity implies
$$\dim \cH_r=\dim (\overline{\iota}\cdot\iota,\rho_r)=n_r.$$
Let $\{V(r)_i\}_{i=1}^{n_r}$ be an orthonormal basis of $\cH_r$.
Then, as in \cite[Section 3]{ILP}, we can show that
$\{(\sqrt{d(r)}V(r)_i^*, \sqrt{d(r)}V(r)_i)\}_{r,i}$
is a quasi-basis for $E$.
For $x\in A$, we set $x(r)_i=d(r)E(V(r)_ix)$.
Then, we have
$$x=\sum_{r\in \cS}\sum_{i=1}^{n_r}V(r)_i^*x(r)_i.$$

Let $C$ be an intermediate $C^*$-subalgebra between $A\supset B$.
Thanks to Theorem 3.3 and the assumption $M(A)\cap B'=\C$, $C$ is simple.
We set
$$\cK_r=\cH_r \cap M(C).$$
We arrange the orthonormal basis $\{V(r)_i\}_{i=1}^{n_r}$ such that
$\{V(r)\}_{i=1}^{m_r}$ is an orthonormal basis of $\cK_r$.

\begin{lemma} Let the notations be as above.
Then, we have
$$C=\{x\in A; \;x(r)_i=0, \; \forall i>m_r,\; r\in \cS\}.$$
\end{lemma}

\begin{proof} Let $D$ be the right-hand side of the above equality.
Then, $C\supset D$ is obvious.
Let $F$ be a conditional expectation from $A$ onto $C$ whose existence
is assured by Proposition 6.1.
We claim that the restriction of $F^M$ to $\cH_r$ is the orthogonal
projection onto $\cK_r$.
Indeed, since $F$ is a $C$-$C$ bimodule map, $F^M$ globally
preserves $\cH_r$.
For $V_1, V_2\in \cH_r$, we have $\inpr {F^M(V_1)}{V_2}=V_2^*F^M(V_1)$.
Since it is already a scalar, it is equal to
$$F^M(V_2^*F^M(V_1))=F^M(V_2^*)F^M(V_1)=F^M(F^M(V_2^*)V_1)
=\inpr {V_1}{F^M(V_2)},$$
which shows the claim.
This implies that for every $x$, we get
$$F(x)=\sum_{r\in \cS}\sum_{i=1}^{n_r}F(V(r)_i^*)x(r)_i
=\sum_{r\in \cS}\sum_{i=1}^{m_r}V(r)_i^*x(r)_i\in D.$$
Therefore, we get the statement.
\end{proof}

The above lemma immediately implies the following:

\begin{corollary} Let $A\supset B$ be an irreducible inclusion of
$C^*$-algebras in $\cC_1$ (or $\cC_2$), and
$E$ be a conditional expectation from $A$ onto $B$ of a finite index.
Then, there exists a one-to-one correspondence between the following two
objects:
$(1)$ The set of intermediate $C^*$-subalgebras.
$(2)$ The set of systems of subspaces $\cK_r\subset \cH_r$, $r\in \cS$
satisfying
\begin{itemize}
\item[$(\mathrm{i})$] $\cK_r \cK_s\subset
\sum_{t\in \cS}B\cK_t$
\item[$(\mathrm{ii})$] $\cK_r^*\subset B\cK_{\overline{r}}$.
\end{itemize}
The correspondence is given as follows: For an intermediate
$C^*$-subalgebra $C$,
$\cK_r$ is given by $\cK_r=M(C)\cap \cH_r$.
On the other hand, for a given system of subspaces $\cK_r\subset \cH_r$,
$r\in \cS$, the corresponding subalgebra is $\sum_{r\in \cS}B\cK_r$.
\end{corollary}

As in \cite{ILP}, we apply this to depth 2 inclusions.
Let $A\supset B$ be an irreducible inclusion of simple $\sigma$-unital
$C^*$-algebras of a finite index, and
$$B\subset A\subset A_1\subset A_2\subset \cdots,$$
be the tower.
We say that $A\supset B$ is of depth 2 if $M(A_2)\cap B'$ is a factor.

\begin{corollary} Let $A\supset B$ be a depth 2 inclusion of simple
$\sigma$-unital $C^*$-algebras.
Then, there exists an action $\delta$ of a finite  dimensional
Kac algebra $Q$ such that $B$ is the fixed point algebra $A^\delta$.
There exists a one-to-one correspondence between the set of intermediate
$C^*$-subalgebras between $A\supset B$ and the set of left coideal
$*$ subalgebras of $Q$.
\end{corollary}

\begin{proof} Let $e$ be a minimal projection of $\K$.
If the statement holds for $(A\otimes \K)\supset (B\otimes \K)$, the
restriction of $\delta$ to $(1\otimes e)(A\otimes \K)(1\otimes e)$
satisfies the conditions of the statement.
Therefore, we may assume that $A$ and $B$ are stable.
The rest of the statement follows from exactly the same arguments
as in \cite[Section 4]{I98} and \cite[Section 4]{ILP}.
\end{proof}

\begin{remark} (1) As is pointed out in \cite[Theorem 5.12]{KajW},
in general, depth two inclusions of $C^*$-algebras cannot be
characterized by crossed products instead of fixed point algebras
as above.
(2) One can get a similar statement for a group-subgroup inclusion
using arguments in \cite{I97}.
\end{remark}

The Galois correspondence of group actions on $C^*$-algebras has been
discussed by several authors \cite{Ch}, \cite{LOP} (see also \cite{NT}
and references in it).
For a finite group action on a simple $C^*$-algebra, we get the following
as a special case of the previous one:

\begin{corollary} Let $\alpha$ be an outer action of a finite group $G$
on a simple $\sigma$-unital $C^*$-algebra $A$.
Then,
\begin{itemize}
\item [$(1)$] There exists a one-to-one correspondence between the set
of subgroups of $G$ and the set of intermediate $C^*$-subalgebras
between $A\rtimes_\alpha G\supset A$.
The correspondence is given as follows: For a subgroup $H$, the
corresponding intermediate $C^*$-subalgebra is $A\rtimes_\alpha H$.
\item [$(2)$] There exists a one-to-one correspondence between the set
of subgroups of $G$ and the set of intermediate $C^*$-subalgebras
between $A\supset A^\alpha$.
The correspondence is given as follows: For a subgroup $H$, the
corresponding intermediate $C^*$-subalgebra is the fixed point subalgebra
of $A$ under the restriction of $\alpha$ to $H$.
\end{itemize}
\end{corollary}

\section{A Kishimoto's type theorem}
Let $B$ be a simple $C^*$-algebra and $\alpha$ be an outer automorphism
of $B$.
In \cite{Ki}, A. Kishimoto proved the following statement
(as a special case):
For every element $x\in M(B)$ and every non-zero hereditary $C^*$-subalgebra
$C$ of $B$,
$$\inf\{||cx\alpha(c)||;\; c\in C_+,\; ||c||=1\}=0.$$
There are a number of applications of this result: for example,
it implies that the reduced crossed product of a simple
(and purely infinite) $C^*$-algebra by an outer action of a
discrete group is again simple (and purely infinite) \cite{Ki}, \cite{KK}.

As observed in \cite{O}, a counterpart of Kishimoto's theorem
for a conditional expectation $E$ from $A$ onto $B$ would be
the following statement:
under some condition corresponding to outerness, for every element
$x\in M(A)$ and every non-zero hereditary $C^*$-subalgebra $C$ of $B$,
$$\inf\{||c(x-E(x))c||;\; c\in C_+,\; ||c||=1\}=0.$$
In fact, Osaka adopted this property as a definition of outerness of
a conditional expectation in \cite{O}.
The goal of this section is to prove it for an irreducible inclusions
of simple $\sigma$-unital $C^*$-algebras of a finite index and finite depth.
In the  course of preparation, we prove Proposition 7.3, which is probably
of enough interest in its own right.

For $A\in \cC_1$, we use the notations $\En A_0=\Mo AA$,
$\mathrm{Sect}(A)=\Se AA$.
$\rho\in \En A_0$ is said to be irreducible if
$M(A)\cap \rho(A)'=\C$ holds.
We say that $\cT=\{\rho_r\}_{r\in \cT_0}\subset \En A_0$ is a closed
system of irreducible endomorphisms if the following properties
are satisfied:
\begin{itemize}
\item [(1)] $[\rho_r]=[\rho_s]$ implies $r=s$.
\item [(2)] There exists $e\in \cT_0$ such that $\rho_e=\Id$.
\item [(3)] There exist non-negative integers $N_{r,s}^t$ for
$r,s,t\in \cT_0$ such that
$$[\rho_r][\rho_s]=\bigoplus_{t\in \cT_0}
N_{r,s}^t[\rho_t].$$
\item [(4)] For each $r\in \cT_0$, there exists
$\overline{r}\in \cT_0$ such that
$\overline{[\rho_r]}=[\rho_{\overline{r}}]$.
\end{itemize}

For each $r\in \cT_0$, we choose
$R_r\in (\Id,\rho_{\overline{r}}\cdot \rho_r)$ and
$\bR_r\in (\Id, \rho_r\cdot\rho_{\overline{r}})$ as in Lemma 4.4, and
for $r,s,t\in \cT_0$, we choose an orthonormal basis
$\{T(_{r,s}^t)_{a}\}_{a=1}^{N_{r,s}^t}$ of $(\rho_t,\rho_r\cdot\rho_s)$.

Now, we assume that $A$ is acting on a Hilbert space $H$ irreducibly,
and set
$$\cL_r=\{V\in \B(H);\; Vx=\rho_r(x)V\}.$$
Thanks to Lemma 5.2, $n_r:=\dim \cL_r$ is finite.
We introduce an inner product into $\cL_r$ by
$$\inpr {V_1}{V_2}1=V_2^*V_1.$$
Lemma 5.1, (1) shows that we can define antilinear maps
$F_r:\cL_r\longrightarrow \cL_{\overline{r}}$ and its inverse
$\overline{F}_r: \cL_{\overline{r}}\longrightarrow \cL_r$ by
$$F_r(V)=\sqrt{d(r)}V^*\bR_r,$$
$$\overline{F}_r(W)=\sqrt{d(r)}W^*R_r.$$
Using these, we can introduce
another inner product into $\cL_r$ by
$$(V_1,V_2)=\inpr{F(V_2)}{F(V_1)}=d(r)\bR_r^*V_1V_2^*\bR_r.$$
We choose an orthogonal basis $\{V(r)_i\}_{i=1}^{n_r}$ of $\cL_r$
with respect to the inner product $\inpr{\cdot}{\cdot}$ that satisfies
$$(V(r)_i,V(r)_j)=c(r)_i\delta_{i,j}.$$

\begin{lemma} Let the notations be as above.
We set
$$C_{(r,i),(s,j)}^{(t,k)}=\sum_{a=1}^{N_{r,s}^t}
\inpr{T(^t_{r,s})_a^*V(r)_iV(s)_j}{V(t)_k}
T(_{r,s}^t)_a\in (\rho_t,\rho_r\cdot\rho_s).$$
Then, the following hold:
\begin{itemize}
\item [$(1)$]
$$V(r)_iV(s)_j=\sum_{t,k}C_{(r,i),(s,j)}^{(t,k)}V(t)_k.$$
\item [$(2)$]
$$V(s)_jV(t)_k^*=\sum_{r,i}\frac{d(r)c(s)_j}{d(t)}V(r)_i^*
C_{(r,i),(s,j)}^{(t,k)}.$$
\end{itemize}
\end{lemma}

\begin{proof} (1). $\sum_{t,a}T(_{r,s}^t)_aT(_{r,s}^t)_a^*=1$ implies
$$V(r)_iV(s)_j=\sum_{t,a}T(_{r,s}^t)_aT(_{r,s}^t)_a^*V(r)_iV(s)_j.$$
Since $T(_{r,s}^t)_a^*V(r)_iV(s)_j\in \cL_t$, we get
$$T(_{r,s}^t)_a^*V(r)_iV(s)_j
=\sum_{t,k}\inpr{T(_{r,s}^t)_a^*V(r)_iV(s)_j}{V(t)_k}V(t)_k,$$
which shows the statement.

(2). Using Lemma 5.1, (1) and a similar argument as above, we get
\begin{eqnarray*}V(s)_jV(t)_k^*
&=&\sqrt{d(s)}F_s(V(s)_j)^*R_sV(t)_k^*
=\sqrt{d(s)}F_s(V(s)_j)^*V(t)_k^*\rho_t(R_s)\\
&=&\sqrt{d(s)}\sum_{r,a}
F_s(V(s)_j)^*V(t)_k^*T(_{t,\overline{s}}^r)_aT(_{t,\overline{s}}^r)_a^*
\rho_t(R_s)\\
&=&\sqrt{d(s)}\sum_{r,i,a}
F_s(V(s)_j)^*V(t)_k^*T(_{t,\overline{s}}^r)_aV(r)_iV(r)_i^*
T(_{t,\overline{s}}^r)_a^*\rho_t(R_s)\\
\end{eqnarray*}
Thanks to the Frobenius reciprocity, we may replace the orthonormal basis
$\{T(_{t\overline{s}}^r)\}_{a=1}^{N_{t\overline{s}}^r}$ with
$$\{\sqrt{\frac{d(r)d(s)}{d(t)}}T(_{r,s}^t)_b^*\rho_r(\bR_s)
\}_{b=1}^{N_{r,r}^t}.$$
Thus,
\begin{eqnarray*}V(s)_jV(t)_k^*
&=&\frac{d(r)\sqrt{d(s)}^3}{d(t)}\sum_{r,i,b}
F_s(V(s)_j)^*V(t)_k^*T(_{r,s}^t)_b^*\rho_r(\bR_s)V(r)_i\\
&\times& V(r)_i^*\rho_r(\bR_s)^*T(_{r,s}^t)_b\rho_t(R_s)\\
&=&\frac{d(r)\sqrt{d(s)}^3}{d(t)}\sum_{r,i,b}
F_s(V(s)_j)^*V(t)_k^*T(_{r,s}^t)_b^*V(r)_i\bR_s\\
&\times& V(r)_i^*\rho_r(\bR_s^*\rho_s(R_s))T(_{r,s}^t)_b\\
&=&\frac{d(r)\sqrt{d(s)}}{d(t)}\sum_{r,i,b}
F_s(V(s)_j)^*V(t)_k^*T(_{r,s}^t)_b^*V(r)_i\bR_sV(r)_i^*T(_{r,s}^t)_b.
\end{eqnarray*}
Since $V(r)_i^*T(_{r,s}^t)_bV(t)_k\in \cL_s$, we have
\begin{eqnarray*}
V(r)_i^*T(_{r,s}^t)_bV(t)_k
&=&\sum_{l}\inpr{V(r)_i^*T(_{r,s}^t)_bV(t)_k}{V(s)_l}V(s)_l\\
&=&\sum_{l}\overline{\inpr{T(_{r,s}^t)_b^*V(r)_iV(s)_l}{V(t)_k}}V(s)_l.
\end{eqnarray*}
Thus,
\begin{eqnarray*}V(s)_jV(t)_k^*
&=&\frac{d(r)\sqrt{d(s)}}{d(t)}\sum_{r,i,l,}
F_s(V(s)_j)^*V(s)_l^*\bR_s
V(r)_i^*C_{(r,i),(s,l)}^{(t,k)}\\
&=&\frac{d(r)}{d(t)}\sum_{r,i,l,}
F_s(V(s)_j)^*F_s(V(s)_l)V(r)_i^*C_{(r,i),(s,l)}^{(t,k)}\\
&=&\frac{d(r)c(s)_j}{d(t)}\sum_{r,i}
V(r)_i^*C_{(r,i),(s,j)}^{(t,k)},
\end{eqnarray*}
which shows the statement.
\end{proof}

\begin{lemma} Let the notations be as above.
Then, the following hold:
\begin{itemize}
\item [$(1)$] Let $r_i,s_j,t\in \cT_0$, $i=1,2,\cdots m$, $j=1,2,\cdots, n$.
If $x(t)_k$ belongs to
$$(\rho_t\cdot\rho_{r_1}\cdots\rho_{r_m},\rho_{s_1}\cdots\rho_{s_n}),$$
for every $t\in\cT_0$, $1\leq k\leq n_k,$ and
$$\sum_{t,k}x(t)_kV(t)_k=0$$
holds, then, $x(t)_k=0$ for all $t,k$.
\item [$(2)$] When $\inpr{T(^t_{r,s})_a^*V(r)_iV(s)_j}{V(t)_k}\neq 0$,
we have $c(r)_ic(s)_j=c(t)_k$.
\item [$(3)$]
$$\sum_{r,i}\frac{d(r)}{c(r)_i}\inpr{C_{(r,i),(s,j)}^{(t,k)}}
{C_{(r,i),(s,j)}^{(t,k')}}
=\frac{d(t)}{c(t)_k}\delta_{k,k'}.$$
\end{itemize}
\end{lemma}

Before starting the proof, we remark an outcome of (1) above.
Assume that we are given $V(s_1)_{j_1}V(s_2)_{j_2}\cdots V(s_n)_{j_n}
V(r_m)_{i_m}^*\cdots V(r_2)_{i_2}^*V(r_1)_{i_1}^*$.
Using Lemma 7.1, (1), (2), and  Lemma 5.1, (1) repeatedly,
we can reduce this element to the form $\sum_{r,i}y(r)_iV(r)_i$,
$y(r)_i\in A$.
Lemma 7.2, (1) asserts that no matter what way we take in the reduction
process, the resulting ``coefficients" $\{y(r)_i\}$ are uniquely determined.

\begin{proof} (1). We show $x(e)=0$.
The general statement can be obtained by applying this and
Lemma 7.1, (2) to
$$\sum_{t,k}x(t)_kV(t)_kV(r)_i^*=0.$$
We choose orthonormal bases
$$\{X(r)_a\}_a\subset (\rho_r,\rho_{r_1}\cdots\rho_{r_m}),$$
$$\{Y(r)_b\}_b\subset (\rho_r,\rho_{s_1}\cdots \rho_{s_n}).$$
Then,
\begin{eqnarray*}
0&=&\bR_r^*Y(r)_b^* \big(\sum_{t,k}x(t)_kV(t)_k\big) X(r)_a\bR_r
=\sum_{t,k}\bR_r^*Y(r)_b^* x(t)_k\rho_t(X(r)_a\bR_r)V(t)_k\\
&=&\bR_r^*Y(r)_b^* x(e) X(r)_a\bR_r,
\end{eqnarray*}
where we use
$$\bR_r^*Y(r)_b^* x(t)_k\rho_t(X(r)_a\bR_r)\in
(\rho_t,\Id)=\delta_{t,e}\C.$$
Since $Y(r)_b^* x(e) X(r)_a\in (\rho_r,\rho_r)$ is already a scalar,
we get $Y(r)_b^* x(e) X(r)_a=0$.
If $r\neq s$, we have $Y(s)_b^* x(e) X(r)_a\in (\rho_r,\rho_s)=\{0\}$.
Therefore,
$$x(e)=\sum_{r,s,a,b}Y(s)_bY(s)_b^*x(e)X(r)_aX(r)_a^*=0.$$

(2). Using Lemma 7.1, (1) and (2), we get
\begin{eqnarray*}
V(r)_iV(s)_jV(t)_k^*
&=&\sum_{t',k'}C_{(r,i),(s,j)}^{(t',k')}V(t')_{k'}V(t)_k^*\\
&=&\sum_{t',r',k',i'}\frac{d(r')c(t')_{k'}}{d(t)}V(r')_{i'}^*
\rho_{r'}(C_{(r,i),(s,j)}^{(t',k')})
C_{(r',i'),(t',k')}^{(t,k)}.
\end{eqnarray*}
On the other hand, using Lemma 7.1, (2) twice, we get
\begin{eqnarray*}
V(r)_iV(s)_jV(t)_k^*
&=&\sum_{u,l}\frac{d(u)c(s)_j}{d(t)}V(r)_iV(u)_l^*
C_{(u,l),(s,j)}^{(t,k)}\\
&=&\sum_{u,r',l,i'}\frac{d(r')c(r)_ic(s)_j}{d(t)}
V(r')_{i'}^*C_{(r',i'),(r,i)}^{(u,l)}
C_{(u,l),(s,j)}^{(t,k)}.
\end{eqnarray*}
Now, we compare the ``$V(e)^*$-coefficient" by applying (1),
and get
$$c(t)_kC_{(r,i),(s,j)}^{(t,k)}=c(r)_ic(s)_jC_{(r,i),(s,j)}^{(t,k)},$$
which shows the statement.

(3). Lemma 7.1, (2) implies
$$V(t)_lV(t)_k^*=\sum_{r,r',i,i'}\frac{d(r)d(r')c(s)_j^2}{d(t)^2}
{C_{(r',i'),(s,j)}^{(t,l)}}^*V(r')_{i'}V(r)_i^*
C_{(r,i),(s,j)}^{(t,k)}.$$
Using Lemma 7.1, (2) again and applying (1) in the same way as above,
we get
\begin{eqnarray*}
\frac{c(t)_k}{d(t)}\delta_{k,l}
&=&\sum_{r,i}\frac{d(r)c(r)_ic(s)_j^2}{d(t)^2}
{C_{(r,i),(s,j)}^{(t,l)}}^*C_{(r,i),(s,j)}^{(t,k)}\\
&=&\sum_{r,i}\frac{d(r)c(t)_k^2}{c(r)_id(t)^2}
{C_{(r,i),(s,j)}^{(t,l)}}^*C_{(r,i),(s,j)}^{(t,k)},\\
\end{eqnarray*}
which finishes the proof.
\end{proof}

\begin{proposition} Let $A$ be a $C^*$-algebra in $\cC_1$ (or $\cC_2$)
acting irreducibly on a Hilbert space $H$, and
$\cT=\{\rho_r\}_{r\in \cT_0}\subset \En A_0$ be a finite
closed system of irreducible endomorphisms.
We set
$$\cL_r=\{V\in \B(H);\; Vx=\rho_r(x)V\}.$$
Let $D$ be the $C^*$-algebra generated by $\cup_{r\in \cT_0}A\cL_r$.
Then, $D$ is simple and there exists a conditional expectation $\cE$
of a finite index from $D$ onto $A$ satisfying
$\cE(aV)=0$ for every $a\in A$, $V\in \cL_r$, $r\neq e$.
\end{proposition}

\begin{proof} Our proof is based on R. Longo's $Q$-system \cite{L94}
and the characterization of the $C^*$-basic construction.
We set $e(r)_i=d(r)/c(r)_i$ and
$$d=\sqrt{\sum_{r,i}e(r)_i}.$$
Since $A$ is stable, we can choose isometries
$\{W(r)_i\}_{r,i}\subset M(A)$ such that
$$\sum_{r,i}W(r)_iW(r)_i^*=1.$$
We define an endomorphism $\gamma$ of $A$ by
$$\gamma(x)=\sum_{r,i}W(r)_i\rho_r(x)W(r)_i^*.$$
Let $S=W(e)$ and
$$T=\frac{1}{d}\sum_{r,s,t,i,j,k}\sqrt{\frac{e(r)_ie(s)_j}{e(t)_k}}
W(r)_i\rho_r(W(s)_j)C_{(r,i),(s,j)}^{(t,k)}W(t)_k^*.$$
Then, thanks to Lemma 7.2, (3),
$S\in (\Id,\gamma)$ and $T\in (\gamma,\gamma^2)$ are isometries
satisfying
$$S^*T=\gamma(S^*)T=\frac{1}{d}.$$
We claim that $TT^*=\gamma(T^*)T$ holds.
Indeed, straightforward calculation yields
\begin{eqnarray*}
TT^*&=&\frac{1}{d^2}\sum_{r,s,t,u,i,j,k,l}
\frac{\sqrt{e(r)_ie(s)_je(u)_le(v)_m}}{e(t)_k}
W(r)_i\rho_r(W(s)_j)C_{(r,i),(s,j)}^{(t,k)}{C_{(u,l),(v,m)}^{(t,k)}}^*\\
&\times& \rho_u(W(v)_m^*)W(u)_l^*.
\end{eqnarray*}
\begin{eqnarray*}
\gamma(T^*)T&=&\frac{1}{d}\sum_{r,t,u,i,k,l}
\sqrt{\frac{e(r)_ie(t)_k}{e(u)_l}}
W(r)_i\rho_r(T^*W(t)_k)C_{(r,i)(t,k)}^{(u,l)}W(u)_l^*\\
&=&\frac{1}{d^2}\sum_{r,s,t,u,v,i,j,k,l,m}
\sqrt{\frac{e(r)_ie(t)_k^2e(v)_m}{e(s)_je(u)_l}}
W(r)_i\rho_r(W(s)_j{C_{(t,k)(v,m)}^{(s,j)}}^*\rho_t(W(v)_m^*))\\
&\times&C_{(r,i)(t,k)}^{(u,l)}W(u)_l^*\\
&=&\frac{1}{d^2}\sum_{r,s,t,u,v,i,j,k,l,m}
\sqrt{\frac{e(r)_ie(t)_k^2e(v)_m}{e(s)_je(u)_l}}
W(r)_i\rho_r(W(s)_j{C_{(t,k)(v,m)}^{(s,j)}}^*)\\
&\times&C_{(r,i)(t,k)}^{(u,l)}\rho_u(W(v)_m^*)W(u)_l^*.\\
\end{eqnarray*}
Therefore, the claim is equivalent to
$$\sum_{t,k}\frac{1}{e(t)_k}
C_{(r,i),(s,j)}^{(t,k)}{C_{(u,l),(v,m)}^{(t,k)}}^*
=\sum_{t,k}
\frac{e(t)_k}{e(s)_je(u)_l}
\rho_r({C_{(t,k)(v,m)}^{(s,j)}}^*)C_{(r,i)(t,k)}^{(u,l)}.$$
This follows from comparison of the $V(e)$-coefficients in two different
ways of expansion of $V(r)_iV(s)_jV(v)_m^*V(u)_l^*$.
Therefore, the claim is proven.
$TT^*=\gamma(T^*)T$ also implies $T^2=\gamma(T)T$ because
$$(T^2-\gamma(T)T)^*(T^2-\gamma(T)T)
=2-T^{2*}\gamma(T)T-T^*\gamma(T^*)T^2=0.$$

We set
$$B=\{x\in A;\; Tx=\gamma(x)T,\; Tx^*=\gamma(x^*)T\}.$$
Then, $B$ is an intermediate $C^*$-subalgebra between
$A\supset \gamma(A)$ with $T\in M(B)$.
For $x\in A$, we define $E(x)=T^*\gamma(x)T$, which is a conditional
expectation from $A$ onto $B$ with a quasi-basis $\{(dS^*,dS)\}$.
We claim that $A\supset B$ is irreducible, and consequently, $B$ is simple.
Indeed, let $\iota$ be the inclusion map from $B$ into $A$.
Then, Lemma 4.4 implies that $[\overline{\iota}][\iota]=[\gamma]$.
Thus, the Frobenius reciprocity implies
$$\dim(\iota,\iota)=\dim(\Id,\gamma)=1,$$
which shows $M(A)\cap B'=\C$.
Theorem 3.3 implies that $B$ is simple.

Now we show $D=A_1$.
Let
$$V=\frac{1}{d}\sum_{r,i}\sqrt{e(r)_i}W(r)_iV(r)_i,$$
which is an isometry satisfying $Vx=\gamma(x)V$ for $x\in A$.
We show that $p=VV^*$ satisfies the condition of Proposition 2.5, (1).
Direct computation using Lemma 7.1, (2) shows
\begin{eqnarray*}
V^*T&=&\frac{1}{d^2}\sum_{r,s,t,i,j,k}\sqrt{\frac{e(r)_i^2e(s)_j}{e(t)_k}}
V(r)_i^*\rho_r(W(s)_j)C_{(r,i),(s,j)}^{(t,k)}W(t)_k^*\\
&=&\frac{1}{d^2}\sum_{r,s,t,i,j,k}\sqrt{\frac{e(r)_i^2e(s)_j}{e(t)_k}}
W(s)_jV(r)_i^*C_{(r,i),(s,j)}^{(t,k)}W(t)_k^*\\
&=&\frac{1}{d^2}\sum_{s,t,j,k}\sqrt{e(t)_ie(s)_j}
W(s)_jV(s)_jV(t)_k^*W(t)_k^*\\
&=&VV^*.
\end{eqnarray*}
On the other hand, Lemma 7.1, (1) implies
\begin{eqnarray*}
TV&=&\frac{1}{d^2}\sum_{r,s,t,i,j,k}\sqrt{e(r)_ie(s)_j}
W(r)_i\rho_r(W(s)_j)C_{(r,i),(s,j)}^{(t,k)}V(t)_k\\
&=&\frac{1}{d^2}\sum_{r,s,i,j}\sqrt{e(r)_ie(s)_j}
W(r)_i\rho_r(W(s)_j)V(r)_iV(s)_j\\
&=&\frac{1}{d^2}\sum_{r,s,i,j}\sqrt{e(r)_ie(s)_j}
W(r)_iV(r)_iV(s)_jV(s)_j\\
&=&V^2,
\end{eqnarray*}
which shows $pT=Tp=V^2V^*$.
Therefore, we get $p\in B'$ because $p$ commutes with $\gamma(A)$ and
$B=T^*\gamma(A)T$.
As we have $A=BS$, to prove $pxp=E(x)p$, it suffices to show that
$pSp$ is equal to
$$E^M(S)p=T^*\gamma(S)Tp=\frac{1}{d}Tp=\frac{1}{d}V^2V^*,$$
which can be shown by direct computation.
Since $B$ is simple, the map $B\ni b\mapsto bp$ is faithful.
Therefore, Proposition 2.5, (1) implies that we can identify
the basic construction $A_1$ with the $C^*$-algebra generated by
$ApA$.
Obviously, we have an inclusion $A_1\supset D$.
On the other hand,
$$AeA\supset AW(r)_i^*eSA=AV(r)_iA=A\rho_r(A)V(r)_i=AV(r)_i,$$
which shows $D=A_1$.

To finish the proof, it suffices to show that
$V(r)_i\in M(D)$ and to put $\cE=E_1$.
However, this follows immediately from Lemma 7.1.
\end{proof}

\begin{remark} In the above, we actually have $c(r)_i=1$.
Indeed, thanks to the local index formula of the minimal conditional
expectation $E_{\rho_r}\cdot \cE$ \cite{H}, \cite{W}, we have
$$\cE(V(r)_iV(r)_i^*)=\frac{1}{d(r)},$$
which implies $c(r)_i=1$.
\end{remark}

\begin{theorem}
Let $A$ be a $C^*$-algebra in $\cC_1$ (or $\cC_2$),
and $\rho\in \En A_0$ be irreducible with $[\rho]\neq [\Id]$.
If $[\rho]$ and $[\br]$ generate only finitely many irreducibles,
then for every $x\in M(A)$ and for every non-zero hereditary
$C^*$-subalgebra $C$ of $A$,
$$\inf\{||cx\rho(c)||;\; c\in C_+,\; ||c||=1\}=0$$
holds (i.e. $\rho$ is properly outer in the sense of \cite{BEK}).
\end{theorem}

\begin{proof} Thanks to Kishimoto's result \cite[Lemma 1.1]{Ki},
when $\rho$ is an automorphism, the statement holds.
Therefore, we may assume that $\rho$ is not an automorphism.
Suppose that the statement does not hold.
Then, in exactly the same way as in \cite[Lemma 1.1]{Ki}, we can show that
there would exist $a\in M(A)$, a non-zero hereditary $C^*$-subalgebra
$C\subset A$, and a positive number $\delta$ such that for every unitary
$u\in C+\C1$ and every pure state
$\varphi\in P(A)$ whose support $s(\varphi)$ belongs to $C^{**}$,
the following holds:
\begin{equation}
||s(\varphi)u^*a\rho^{**}(us(\varphi))||\geq \delta
\end{equation}
This implies $\pi_\varphi(\rho^{**}(s(\varphi)))\neq 0$.
Let
$$\Lr=\{V;\; V\pi_\varphi(x)=\pi_\varphi(\rho(x))V,\; x\in A\}.$$
We claim that the map $\Lr\ni V\mapsto V\Omega_\varphi\in H$
gives an isometry from $\Lr$ onto
$\pi_\varphi(\rho^{**}(s(\varphi)))H_\varphi$.
Indeed, let $\psi$ be the pure state on $\rho(A)$ given by
$\varphi\cdot\rho^{-1}$.
Then, $\pi_\varphi(\rho^{**}(s(\varphi)))H_\varphi$ is characterized
as the set of vectors in $H_\varphi$ inducing a scalar multiple of
$\psi$ on $\rho(A)$, which shows that the above map is an isometry
$\Lr$ into $\pi_\varphi(\rho^{**}(s(\varphi)))H_\varphi$.
On the other hand, for a given
$\xi\in \pi_\varphi(\rho^{**}(s(\varphi)))H_\varphi$, we can construct
$V\in \Lr$ by setting
$$V\pi_\varphi(x)\Omega_\varphi=\pi_\varphi(\rho(x))\xi,\quad x\in A.$$
Thus, the claim holds.

Lemma 5.1 and Lemma 5.2 imply that $\cL_\rho$
(and $\pi_\varphi(\rho^{**}(s(\varphi)))H_\varphi$ as well)
is finite dimensional.
Let $\{V_i\}_{i=1}^m$ be an orthonormal basis of $\Lr$.
Then, thanks to the claim, we get
$$\rho^{**}(s(\varphi))=\sum_{i=1}^mV_is(\varphi)V_i^*,$$
which implies
$$s(\varphi)u^*a\rho^{**}(us(\varphi)u^*)a^*us(\varphi)
=\sum_{i=1}^m|\inpr {\pi_\varphi(a)V_i\pi_\varphi(u)\Omega_\varphi}
{\pi_\varphi(u)\Omega_\varphi}|^2 s(\varphi).$$
Let $K$ be the closure of $\pi_\varphi(C)\Omega_\varphi$.
Then, (7.1) is equivalent to
$$\sum_{i=1}^m|\inpr {\pi_\varphi(a)V_i\xi}{\xi}|^2\geq
\delta^2||\xi||^4,\quad \xi\in K.$$

We introduce a map
$\Phi:K\times K\longrightarrow \C^m$ by
$$\Phi(\xi,\eta)=(\inpr {\pi_\varphi(a)V_1\xi}{\eta},
\inpr {\pi_\varphi(a)V_2\xi}{\eta},\cdots,
\inpr {\pi_\varphi(a)V_m\xi}{\eta}),$$
and set
$$\delta_0=\inf\{||\Phi(\xi,\xi)||;\ \xi\in K,\; ||\xi||=1\}\geq \delta.$$
We take a sequence of unit vectors $\{\xi_n\}\in K$ satisfying
$$\lim_{n\to\infty} ||\Phi(\xi_n,\xi_n)||=\delta_0.$$
We fix a free ultrafilter $\omega\in \beta\N\setminus \N$.
For a unit vector $\eta\in K$, we set
$$e=\lim_{n\to \omega}\Phi(\xi_n,\xi_n),$$
$$f=\lim_{n\to \omega}(\Phi(\xi_n,\eta)+\Phi(\eta,\xi_n)),$$
$$g=\Phi(\eta,\eta),$$
$$h=\lim_{n\to \omega}(\inpr{\xi_n}{\eta}+\inpr{\eta}{\xi_n}).$$
Note that since $\Phi$ is bounded on the unit ball of $K$, the above
limits exist.
For all $t\in \R$, we have
\begin{eqnarray*}\lefteqn{
\lim_{n\to\omega}(||\Psi(\xi_n+t\eta,\xi_n+t\eta)||^2-
\delta_0^2||\xi_n+t\eta||^4)}\\
&=&||e+tf+t^2g||^2-\delta_0^2(1+th+t^2)^2\\
&=&(||g||^2-\delta_0^2)t^4
+\big(\inpr{f}{g}+\inpr{g}{f}-\delta_0^2(h+\overline{h})\big)t^3\\
&+&\big(\inpr{e}{g}+\inpr{g}{e}+||f||^2-\delta_0^2(2+|h|^2)\big)t^2
+\big(\inpr{e}{f}+\inpr{f}{e}-\delta_0^2(h+\overline{h})\big)t\\
&\geq&0.
\end{eqnarray*}
This implies
$$\inpr{e}{f}+\inpr{f}{e}=\delta_0^2(h+\overline{h}),$$
$$\inpr{e}{g}+\inpr{g}{e}+||f||^2\geq\delta_0^2(2+|h|^2)\big. $$
>From the latter, we get
$$\inpr{e}{g}+\inpr{g}{e}+2\lim_{n\to\omega}
(||\Phi(\xi_n,\eta)||^2+||\Phi(\eta,\xi_n)||^2)
\geq2\delta_0^2\big. $$
Let $V=\sum_{i=1}^m\overline{e_i}V_i$, and $p_{\xi_n}$ be
the projection onto $\C\xi_n$.
The above inequality is equivalent to
$$\inpr{(\pi_{\varphi}(a)V+V^*\pi_\varphi(a)^*)\eta}{\eta}
+2\inpr{T\eta}{\eta}\geq 2\delta_0^2,$$
where
$$\mathrm{w-}\lim_{n\to\omega}\sum_{i=1}^m
(\pi_\varphi(a)V_ip_{\xi_n}V_i^*\pi_\varphi(a)^*
+V_i^*\pi_\varphi(a)^*p_{\xi_n}\pi_\varphi(a)V_i).$$
Since $C$ is a hereditary $C^*$-subalgebra with $s(\varphi)\in C^{**}$,
$K$ is actually the closure of $\pi_\varphi(C)H_\varphi$.
Therefore, for any $x\in C$ we get
$$\pi_\varphi(x)^*(\pi_{\varphi}(a)V+V^*\pi_\varphi(a)^*+2T)\pi_\varphi(x)
\geq \delta_0^2\pi_\varphi(x^*x).$$

To obtain contradiction applying Proposition 7.3, we need to get
rid of $T$.
For any orthonormal system $\{\zeta_j\}_{j=1}^k$ in $H_\varphi$, we have
\begin{eqnarray*}\lefteqn{
\sum_{j=1}^k\sum_{i=1}^m
\inpr{(\pi_\varphi(a)V_ip_{\xi_n}V_i^*\pi_\varphi(a)^*
+V_i^*\pi_\varphi(a)^*p_{\xi_n}\pi_\varphi(a)V_i)\zeta_j}{\zeta_j}}\\
&=&\sum_{j=1}^k\sum_{i=1}^m(|\inpr{\varphi(a)V_i\xi_n}{\zeta_j}|^2
+|\inpr{V_i^*\pi_\varphi(a)^*\xi_n}{\zeta_j}|^2)\\
&\leq&\sum_{i=1}^m(||\pi_\varphi(a)V_i\xi_n||^2
+||V_i^*\pi_\varphi(a)^*\xi_n||^2)\\
&\leq& 2m||a||^2,
\end{eqnarray*}
which implies that $T$ is a trace class operator.
Let $\pi$ be the quotient map from $\B(H_\varphi)$ onto the Calkin algebra
$\B(H_\varphi)/\K(H_\varphi)$.
Then, we get
$$\pi(\pi_\varphi(x)^*(\pi_{\varphi}(a)V
+V^*\pi_\varphi(a)^*)\pi_\varphi(x))
\geq \delta_0^2\pi(\pi_\varphi(x^*x)).$$

Let $\cT$ be the finite closed system of irreducible endomorphisms in
$\En A_0$ generated by $\rho$ and $\br$, and  $D$ and $\cE$ be as in
Proposition 7.3, where we identify $A$ with its image in $\B(H_\varphi)$.
Then, thanks to Proposition 7.3, $D$ is simple, and so either
$D\cap \K(H_\varphi)$ is trivial or $D=\K(H_\varphi)$.
The latter does not occur because it would imply that
$\pi_\varphi(A)$ would be a proper $C^*$-subalgebra of $\K(H_\varphi)$
acting irreducibly on $H_\varphi$.
Therefore, the restriction of $\pi$ to $D$ is faithful, and so we get
$$\pi_\varphi(x)^*(\pi_{\varphi}(a)V
+V^*\pi_\varphi(a)^*)\pi_\varphi(x)
\geq \delta_0^2\pi_\varphi(x^*x).$$
Applying $\cE$ to the both sides, we obtain
$$0\geq \delta_0^2\pi_\varphi(x^*x),$$
which is contradiction.
Therefore, we finally get the statement.
\end{proof}

\begin{corollary} Let $A\supset B$ be an irreducible and finite depth
inclusion of simple $\sigma$-unital $C^*$-algebras, and $E$ be a
conditional expectation from $A$ onto $B$ of a finite index.
Then, for every element $x\in M(A)$ and every non-zero hereditary
$C^*$-subalgebra $C$ of $B$,
$$\inf\{||c(x-E(x))c||;\; c\in C_+,\; ||c||=1\}=0.$$
\end{corollary}

\begin{proof} It suffices to show the statement for
$(A\otimes \K, B\otimes \K, E\otimes \Id)$ as usual and we may assume
that $A$ and $B$ are stable from the beginning.
Then, as in the usual finite group crossed product case,
the statement follows from repeated use of Theorem 7.5 and the
``crossed product expression" discussed in Section 6.
\end{proof}

\begin{remark}
It might be an interesting problem to seek for a statement corresponding
to \cite[Theorem 1]{BEK} in our situation.
\end{remark}

\end{document}